\documentclass[11pt,reqno,twoside,a4paper]{amsart}

\usepackage{amsmath, amssymb, amscd, amsthm}

\theoremstyle{plain}
\newtheorem{thm}{Theorem}[section]

\newtheorem{lem}[thm]{Lemma}
\newtheorem{prop}[thm]{Proposition}
\newtheorem{Def}[thm]{Definition}

\theoremstyle{remark}
\newtheorem{rem}[thm]{Remark}
\numberwithin{equation}{section}

\newcommand{\R}{\mathbb R}
\newcommand{\N}{\mathbb N}
\newcommand{\C}{\mathbb C}
\newcommand{\Z}{\mathbb Z}
\newcommand{\al}{\alpha}
\newcommand{\be}{\beta}
\newcommand{\ga}{\gamma}
\newcommand{\Ga}{\Gamma}
\newcommand{\de}{\delta}
\newcommand{\De}{\Delta}
\newcommand{\eps}{\varepsilon}
\newcommand{\si}{\sigma}

\newcommand{\Te}{\Theta}
\newcommand{\la}{\lambda}

\newcommand{\hf}{\frac{1}{2}}
\newcommand{\F}[5]{\,_{#1}F_{#2} \left( \genfrac{.}{.}{0pt}{}{#3}{#4}\ ;#5 \right)}
\newcommand{\Res}[1]{\underset{#1}{\mathrm{Res}}\,}

\begin{document}
\title{Fourier transforms related to a root system of rank 1}
\date{\today}
\author{Wolter Groenevelt}
\address{Korteweg-De Vries Instituut voor Wiskunde\\
Universiteit van Amsterdam\\
Plantage Muidergracht 24\\
1018 TV Amsterdam\\
The Netherlands}
\email{wgroenev@science.uva.nl}
\thanks{The author is supported by the Netherlands Organization for Scientific Research (NWO) for the Vidi-project ``Symmetry and modularity in exactly solvable models.'' Parts of this research were done during the authors stay at the Department of Mathematics at Chalmers University of Technology and G\"oteborg University in Sweden, supported by a NWO-TALENT stipendium of the Netherlands Organization for Scientific Research (NWO)}

\begin{abstract}
We introduce an algebra $\mathcal H$ consisting of difference-reflection operators and multiplication operators that can be considered as a $q=1$ analogue of Sahi's double affine Hecke algebra related to the affine root system of type $(C^\vee_1, C_1)$. We study eigenfunctions of a Dunkl-Cherednik-type operator in the algebra $\mathcal H$, and the corresponding Fourier transforms. These eigenfunctions are non-symmetric versions of the Wilson polynomials and the Wilson functions. 
\end{abstract}

\maketitle

\section{Introduction}
Cherednik's double affine Hecke algebras and their degenerate versions are very useful for studying Macdonald orthogonal polynomials and generalized Fourier transforms, see e.g.~Cherednik \cite{Ch95}, \cite{Ch97}, Macdonald \cite{Mac03},  Opdam \cite{Op95}, Sahi \cite{Sa99}, Stokman \cite{St03}. In rank 1 this approach leads to new interpretations and new results for many well-known orthogonal polynomials of (basic) hypergeometric type and corresponding integral transforms. In this paper we are interested in generalized Fourier transforms related to the rank 1 root system of type $(C^\vee, C)$. The kernels in our Fourier transforms are non-symmetric versions of the Wilson polynomials \cite{Wi80} and the Wilson functions \cite{Gr03}.

Sahi \cite{Sa99} defined a double affine Hecke algebra $\mathcal H_q$ related to the non-reduced root system of type $(C^\vee, C)$, using Noumi's  representation of the affine Hecke algebra associated to root systems of type $C$ \cite{No95}. Sahi's algebra in rank 1 is very useful for studying Askey-Wilson polynomials \cite{AW85} and Askey-Wilson functions \cite{KS01}, see Noumi and Stokman \cite{NS04} and Stokman \cite{St03}. In fact, Sahi used his double affine Hecke algebra to study the Koornwinder-Macdonald polynomials, the multivariable versions of the Askey-Wilson polynomials introduced by Koornwinder \cite{Ko92}. In rank 1 the double affine Hecke algebra $\mathcal H_q$ depends, besides $q$, on four independent parameters, corresponding to the number of $\mathcal W$-orbits in the root system, where $\mathcal W$ is the affine Weyl group. These parameters correspond to the four parameters of the Askey-Wilson polynomials. A certain isomorphism of the double affine Hecke algebra leads to the duality property for the Askey-Wilson polynomials and functions. This property basically says that the geometric and spectral variable can be interchanged. In the limit $q\rightarrow 1$ the double affine Hecke algebra goes over into the degenerate double affine Hecke algebra, see \cite{Ch97A} for these limits in case of reduced root systems. This limit transition between algebras corresponds to the limit transition between Askey-Wilson polynomials and Jacobi polynomials. However, the polynomials and the algebra obtained in the limit depend only on two independent parameters instead of four, and there is no longer a duality property. For $q \rightarrow 1$ the Askey-Wilson polynomials also have a limit that depends on all four parameters, namely the Wilson polynomials. It is a natural question to ask if there is a degenerate double affine Hecke algebra related to the Wilson polynomials. 

In this paper we construct such an algebra related to the Wilson polynomials. We define a representation of the affine Weyl group which resembles Noumi's representation of the affine Hecke algebra associated to root systems of type $C$. For this representation we use difference operators, similar to the difference operators used by Cherednik \cite{Ch97A} to obtain the inverse of the Harish-Chandra transform. With our representation of the affine Weyl group, we construct an algebra $\mathcal H$ which is a new example of a degenerate Hecke algebra. It is shown that $\mathcal H$ is a degeneration of Sahi's double affine Hecke algebra of type $(C^\vee, C)$ in rank 1. In particular, $\mathcal H$ depends on four independent parameters. We remark that there exists a degeneration of Sahi's algebra in rank $n$ similar to the algebra $\mathcal H$ studied in this paper. In the limit the $C_n$ braid relations are deformed, so in rank $n$ the algebra $\mathcal H$ is not obtained from the representation of the Weyl group for the root system of type $(C^\vee_n, C_n)$. The rank $n$ case will be studied in a future paper. 

As in Sahi's algebra, we have a duality isomorphism $\si$ for the algebra $\mathcal H$. To this isomorphism $\si$ one can associate a Fourier transform $\mathbb F:V \rightarrow W$, which is a linear map that intertwines the actions of $\mathcal H$ and $\si(\mathcal H)$ on the vector spaces $V$ and $W$ respectively, that is, $\mathbb F \circ X = \si(X) \circ \mathbb F$ for all $X \in \mathcal H$. We consider Fourier transforms that can be written as integral transforms with kernels the eigenfunctions of the Dunkl-Cherednik-type operator. The eigenfunctions of a Dunkl-Cherednik-type operator in our representation are non-symmetric versions of the Wilson polynomials and the Wilson functions. 

We study the non-symmetric Wilson polynomials using the degenerate affine Hecke algebra $\mathcal H$, which is very similar to the study of non-symmetric Askey-Wilson polynomials in \cite{NS04}. The non-symmetric Wilson functions, however, cannot be treated in the same way as the non-symmetric Askey-Wilson functions (the rank 1 Cherednik kernel in \cite{St03}). For instance, the non-symmetric Askey-Wilson functions are defined in \cite{St03} as a sum of non-symmetric Askey-Wilson polynomials, but a similar sum of non-symmetric Wilson polynomials does not converge absolutely.

Recently, Zhang \cite{Zh05} gave an interpretion of the (symmetric) multivariable Wilson polynomials in the representation theory of the degenerate double affine Hecke algebra by showing that the multi-variable Wilson polynomials can be obtained from symmetric multivariable Jacobi polynomials using the Harish-Chandra transform. For rank 1 this was proved by Koornwinder \cite{Ko85}. The methods in this paper to study the Wilson polynomials are different from the method used in \cite{Zh05}. 

The structure of this paper is as follows. In section \ref{sec2} we define the affine Weyl group, construct a representation, and define a degenerate double affine Hecke algebra using this representation. The non-symmetric Wilson polynomials are needed here to decompose our representation into irreducible representations. In section \ref{sec3} we study the non-symmetric Wilson polynomials. We find orthogonality relations with respect to a complex measure, we obtain the duality property, and we relate the non-symmetric Wilson polynomials to the well-known $_4F_3$-polynomials defined by Wilson \cite{Wi80}. In section \ref{sec4} we define the polynomial Fourier transform as an integral transform with the non-symmetric Wilson polynomials as a kernel. This Fourier transform maps polynomials to functions that are finitely supported on the polynomial spectrum $\Ga$ of the Dunkl-Cherednik-type operator. The algebra $\mathcal H$ also acts on these functions, and we show that the Fourier transform intertwines actions of the algebra $\mathcal H$ on polynomials with the actions on functions on $\Ga$. We also determine the inverse transform, and we give Plancherel-type formulas. Finally, in section \ref{sec5} we define and study a Fourier transform with a non-polynomial kernel. We first introduce a Gaussian-type function. With this Gaussian, the polynomial Fourier transform and its inverse, we construct a Fourier transform that maps polynomials multiplied by the Gaussian to polynomials multiplied by another Gaussian. This Fourier transform also intertwines actions of the algebra $\mathcal H$. We show that we can write this Fourier transform as an integral transform with a non-symmetric Wilson function as a kernel, by relating it to the Wilson function transform I. This non-polynomial Fourier transform turns out to be self-dual. \\

\emph{Acknowledgement.} I would like to thank Jasper Stokman and Genkai Zhang for helpful suggestions and comments.

\section{A degenerate double affine Hecke algebra} \label{sec2}
In the section we define a representation of the affine Weyl group corresponding to the affine root system $(C^\vee_1, C_1)$. With this representation we construct a degenerate double affine Hecke algebra.

\subsection{The affine root system}
Let $E$ be the space of affine linear transformation of $\R$. 
Let $\de$ be the constant function $1$, then for every $f \in E$ we may write
\[
f(u) = vu+\mu \de(u), \qquad u \in \R,
\]
for some $v,\mu \in \R$, so we can identify $E$ with $\R \oplus \R\de$. 
Let $\mathcal R \subset E$ be the subset 
\[
\mathcal R = \{ \pm 1 + \hf m \de, \pm 2 + m \de \ | \ m \in \Z\}.
\]
The set $\mathcal R$ is the non-reduced affine root system of type $(C^\vee_1, C_1)$. 
Define a bilinear form $(\cdot,\cdot)$ on $E$ by 
\[
(v_1+\mu_1 \de, v_2 + \mu_2 \de) = v_1v_2.
\]
For every root $\al \in \mathcal R$ we define $s_\al:E \rightarrow E$ to be the reflection in the affine hyperplane $H_\al=\al^{-1}(0)$, then $s_\al$ is given by 
\[
s_\al(f) = f - (f,\al^\vee)\al, \qquad f \in E,
\]
where $\al^\vee=2\al/(\al,\al)$. The affine Weyl group $\mathcal W$ is the group generated by the reflection $s_\al$, $\al \in \mathcal R$. 

Let $S=\{\pm 2\} \subset \mathcal R$ and $S^\vee=\{\pm 1\} \subset \mathcal R$. The subset $R=S \cup S^\vee \subset \mathcal R$ is a non-reduced (non-affine) root system with Weyl group $W$ generated by the reflections $s_\al, \al \in R$. We set $a_1=2 \in \mathcal R$ and $a_0=\de-2 \in \mathcal R$, then $a_i^\vee=a_i/2$, $i=0,1$. Then the affine Weyl group $\mathcal W$ is the Coxeter group with two generators $s_0=s_{a_0}$, $s_1=s_{a_1}$ and relations $s_i^2=1$. The subgroup $W$ is generated by $s_1$, so $W=\{1,s_1\}$. Observe that $s_1s_0 \in \mathcal W$ acts as a translation operator
\[
(s_1s_0 f)(u)=f(u+1), \qquad f \in E,\ u \in \R.
\]
Setting $\tau(\la) = (s_1 s_0)^\la \in \mathcal W$, for $\la \in \Z$, we obtain a description of $\mathcal W$ as the semi-direct product
\[
\mathcal W = W\ltimes \tau(\Z).
\]

The affine root system $\mathcal R$ decomposes into four $\mathcal W$-orbits, namely
\begin{gather*}
\mathcal Wa_0 = S+(1+2\Z)\de,  \qquad \mathcal Wa_1 = S +2\Z\de,\\
\mathcal Wa_0^\vee = S^\vee+(\hf+\Z)\de,  \qquad \mathcal Wa_1^\vee = S^\vee +\Z\de.
\end{gather*}
Let $\mathcal S=\mathcal Wa_0 \cup \mathcal Wa_1$ and $\mathcal S^\vee =  \mathcal Wa_0^\vee \cup \mathcal Wa_1^\vee$. Then $\mathcal S$ and $\mathcal S^\vee$ are reduced affine root systems with basis $\{a_0,a_1\}$, respectively $\{a_0^\vee, a_1^\vee\}$, and $\mathcal R$ is the disjoint union $\mathcal S \cup \mathcal S^\vee$. Both $\mathcal S$ and $\mathcal S^\vee$ have $\mathcal W$ as corresponding affine Weyl group.

\subsection{The group algebra $\C[\mathcal W]$}
We extend the action of $\mathcal W$ on $E$ to an action of $\mathcal W$ on meromorphic functions on $\C$ by
\[
(s_0 f)(x) = f(1-x), \qquad
(s_1 f)(x) = f(-x).
\]
In particular, this gives an action of $\mathcal W$ on the algebra $\mathcal A$ of polynomials in one variable. Using the expression $\tau(1)= s_1 s_{0}$ it follows that $\tau(\la)$, $\la \in \Z$, acts as the translation operator
\[
(\tau(\la)f)(x) = f( x+\la).
\]

Let $\mathbf t$ be a multiplicity function, i.e., a function $\mathbf t: \mathcal R \rightarrow \C$ that is constant on $\mathcal W$-orbits. Since $\mathbf t$ is completely determined by the values $t_{a_0}$, $t_{a_1}$, $t_{a_0^\vee}$ and $t_{a_1^\vee}$, we may identify $\mathbf t$ with the ordered $4$-tuple $(t_{a_0}, t_{a_0^\vee}, t_{a_1}, t_{a_1^\vee})$. Moreover, to ease notations we will write
\[
\begin{split}
t_{a_0} = t_0, \quad t_{a_0^\vee}=u_0, \quad
t_{a_1} = t_1, \quad t_{a_1^\vee}=u_1.
\end{split}
\]
For an object $O(\mathbf t)$ (function, operator, etc.) depending on a multiplicity function $\mathbf t$ we will often suppress the dependence on $\mathbf t$ and we write $O(\mathbf t)=O$.

Let $\C[\mathcal W]$ be the group algebra of the affine Weyl group $\mathcal W$.
We assume $t_0,t_1 \neq 0$. In $\C[\mathcal W]$ we define
\[
\si_0=t_0 s_0, \qquad \si_1=t_1 s_1, \qquad \xi= \si_0+\si_1.
\]
From $ (\xi-\si_i)^2 = t_{1-i}^2$ it follows that
\begin{equation} \label{eq:xi squared}
\xi^2 = \si_i\xi +\xi \si_i + t_{1-i}^2 - t_{i}^2, \qquad i=0,1. 
\end{equation}
Now we see that $\xi^2$ commutes with $\si_0$ and $\si_1$. Moreover, iterating \eqref{eq:xi squared} for $i=1$ gives the relations
\[
\si_1\xi^{2n}-\xi^{2n}\si_1 =0, \quad \si_1 \xi^{2n+1} + \xi^{2n+1}\si_1 = (\xi^2-t_0^2+t_1^2)\xi^{2n},\qquad n \in \N.
\]
Let $\mathcal A_\xi$ denote the commutative subalgebra of $\C[\mathcal W]$ consisting of elements $p(\xi)$, for $p \in \mathcal A$. Then the above relations can be written as
\begin{equation} \label{eq:S1 xi} 
\si_1 p(\xi) -(s_1p)(\xi) \si_1 = \frac{ t_1^2-t_0^2+\xi^2 }{2\xi}\Big(p(\xi)-(s_1p)(\xi)\Big),
\end{equation}
where $\frac1\xi=\xi^{-1}$. Note that the right hand side is an element of $\mathcal A_\xi$. This immediately leads to the following proposition.
\begin{prop} \label{prop:iso}
The sets $\{\xi^m, \xi^n \si_1\ | \ n,m \in \Z_{\geq 0}\}$ and $ \{\xi^m,\si_1 \xi^n \ | \ n,m \in \Z_{\geq 0}\}$ are linear bases for $\C[\mathcal W]$. So
\[
\mathcal A_\xi \otimes \C[W] \cong \C[\mathcal W] \cong \C[W] \otimes \mathcal A_\xi
\]
as vector spaces. 
\end{prop}
Let $\mathcal A^W\subset \mathcal A$ be the subalgebra of $W$-invariant polynomials, i.e.~$\mathcal A^W=\C[x^2]$, and let $\mathcal A_\xi^W \subset \mathcal A_\xi$ be the subalgebra of $\C[\mathcal W]$ consisting of elements $p(\xi)$ for $p \in \mathcal A^W$. The subalgebra $\mathcal A_\xi^W$ is the center of $\C[\mathcal W]$. This can for example be checked using Proposition \ref{prop:iso} and \eqref{eq:S1 xi}.

\subsection{A representation of $\C[\mathcal W]$} 
To the simple roots $a_0,a_1 \in \mathcal S$ we associate linear operators
$T_{a_i}(\mathbf t)=T_i$ by
\[
\begin{split}
(T_if)(x)& = t_i f(x)+ c_i(x)\Big((s_if)(x)-f(x)\Big), 
\end{split}
\]
where $c_i$ are rational functions defined by
\[
c_i(x;\mathbf t)= \frac{ \big(t_i + u_{i} +a_i^\vee(x)\big)\big(t_i-u_{i} + a_i^\vee(x)\big)} {a_i(x)}.
\]
Recall here that $a_0(x)=1-2x$ and $a_1(x)=2x$. It is easy to verify that the functions $c_i$ satisfy 
\begin{equation} \label{eq:ci}
c_i(x)+(s_ic_i)(x)=2t_i, \qquad i=0,1.
\end{equation}
\begin{prop} \label{prop:Ti}
For $i=0,1$ the operators $T_i$ map $\mathcal A$ into itself.
\end{prop}
To prove the proposition we introduce difference-reflection operators $D_i$ related to the simple roots $a_i \in \mathcal S$,
\[
D_i = \frac{s_i -1}{a_i}.
\]
Proposition \ref{prop:Ti} then follows immediately from the following lemma.
\begin{lem} \label{lem:D}
The operators $D_i$, $i=0,1$, map $\mathcal A$ into itself. Moreover, we have
\[
\begin{split}
D_0 x^m& = \begin{cases}
\hf m x^{m-2} + \mathrm{l.o.t.}& m \text{ even},\\
x^{m-1}-\hf(m-1)x^{m-2}+\mathrm{l.o.t.}& m \text{ odd},
\end{cases}\\
D_1 x^m& =   \begin{cases}
0, & m \text{ even},\\
-x^{m-1}, & m \text{ odd}.
\end{cases}
\end{split}
\]
\end{lem}
Here `l.o.t.' means `lower order terms with respect to the degree'. 
\begin{proof}
By linearity it is enough the prove that $D_ip \in \mathcal A$ for $p(x) = x^m$, $m \in \Z_{\geq 0}$. For $i=1$ the calculation is straightforward. 
For $i=0$ we have
\[
D_0x^m =  \frac{ (1-x)^{m}-x^{m} }{1-2x}.
\]
We replace $x$ by $z+\hf$, then we see that the numerator of the expression
\[
\frac{ (\hf-z)^{m}- (\hf+z)^{m} }{-2z}
\]
only contains odd powers of $z$, so we see that $D_0x^m \in \mathcal A$. To get the explicit expression we expand $ (\hf+z)^m =\sum_{i=0}^m \binom{m}{i}z^i 2^{i-m}$ and write out the above expression.  
\end{proof}

\begin{prop}
The assignment 
\[
\si_i \mapsto T_i, \qquad i=0,1,
\]
extends uniquely to an algebra homomorphism $\pi_{\mathbf t}:\C[\mathcal W] \rightarrow \mathrm{End}(\mathcal A)$.
\end{prop}
\begin{proof}
We only need to verify that $T_i$ satisfies the relations $T_i^2=t_i^2$. Let $p \in \mathcal A$, then we have for $i=0,1$,
\[
\begin{split}
(T_i^2p)(x) =& T_i\Big(t_ip(x) + c_i(x)\big((s_ip)(x) - p(x)\big) \Big)\\
=& t_i\Big(t_ip(x) + c_i(x)\big((s_ip)(x) - p(x)\big)\Big) \\
&+ c_i(x)\Big(t_i(s_ip)(x) + (s_ic_i)(x)\big(p(x) - (s_ip)(x)\big)\Big)\\
&- c_i(x)\Big(t_ip(x) + c_i(x)\big((s_ip)(x) - p(x)\big)\Big)\\
=& \Big(t_i^2 + c_i(x)[c_i(x)+(s_ic_i)(x) -2t_i]\Big)p(x) \\
& - c_i(x)[c_i(x)+(s_ic_i)(x) -2t_i](s_ip)(x).
\end{split}
\]
By \eqref{eq:ci} the terms between square brackets are equal to zero, hence $T_i^2=t_i^2$ as required.
\end{proof}
We denote $\mathcal T=\pi_{\mathbf t}(\C[\mathcal W]) \subset \mathrm{End}(\mathcal A)$, so $\mathcal T$ is the algebra generated by $T_0$ and $T_1$. Moreover, we will denote $T=\pi_{\mathbf t}(\C[W])$, so $T$ is the algebra generated by $T_1$. In the study of the representation $\pi_{\mathbf t}$ the image of $\xi \in \C[\mathcal W]$ will play an important role. We denote 
\[
Y = \pi_{\mathbf t}(\xi) = T_0+T_1 \in \mathcal T.
\]
For $Y$ we have the following useful property. 
\begin{prop} \label{prop:trian}
$Y$ is a triangular operator, i.e., for $m\in \Z_{\geq 0}$,
\[
Y x^m = \ga_m x^m + \mathrm{l.o.t.}
\]
The coefficients $\ga_m$ are given explicitly by
\[
\ga_m = \begin{cases}
t_0+t_1+n, & m=2n,\\
-(t_0+t_1+n), & m=2n-1. 
\end{cases} 
\]
\end{prop}
From here on we assume that $t_0+t_1 \not\in -\hf \N$, so that $\ga_{m_1} \neq \ga_{m_2}$ if $m_1\neq m_2$.
\begin{proof}
We write
\[
\begin{split}
T_0& = t_0 + (t_0+u_0+\hf-x)(t_0-u_0+\hf-x)D_0,\\
T_1& = t_1+ (t_1+u_1+x)(t_1-u_1+x)D_1,
\end{split}
\]
then we obtain from Lemma \ref{lem:D} 
\[
\begin{split}
T_0 x^m& = \begin{cases}
 x^{m+1} - (t_0+\hf(m+1))x^m+\mathrm{l.o.t.}, & m \text{ odd},\\
(t_0+\hf m)x^m+ \mathrm{l.o.t.},& m \text{ even},
\end{cases}\\
T_1 x^m& = \begin{cases}
-x^{m+1}-t_1x^m+\mathrm{l.o.t.}, & m \text{ odd},\\
t_1 x^m +\mathrm{l.o.t.},& m \text{ even}.
\end{cases}
\end{split}
\]
Now the result follows from $Y=T_0+T_1$.
\end{proof}
From Proposition \ref{prop:trian} it follows that for all $m \in \Z_{\geq 0}$ there exists, up to a multiplicative constant, a unique element in $\mathcal A$ which is an eigenfunction of $Y$ for eigenvalue $\ga_m$. We call this element a non-symmetric Wilson polynomial.
\begin{Def} \label{def:Wilson pol}
The non-symmetric Wilson polynomial $p_m$ is the unique monic polynomial of degree $m$ such that $Yp_m = \ga_m p_m$.
\end{Def}
Clearly the set $\{p_m\ | \ m \in \Z_{\geq 0}\}$ forms a linear basis for $\mathcal A$. In order to decompose $\mathcal A$ as an $\mathcal T$-module, we need to know the action of $T_1$ on the non-symmetric Wilson polynomials. By $T_0=Y-T_1$ we then also have the action of $T_0$.
\begin{prop} \label{prop:V1}
We have $T_1p_0=t_1 p_0$, and for $m \geq 1$
\[
T_1 p_m = \begin{cases}
\displaystyle -p_{m+1} + b_m\,p_{m},& m \text{ odd},\\
\displaystyle b_m\, p_m + (b_m^2-t_1^2)\,p_{m-1}, & m \text{ even},
\end{cases}
\]
where
\[
b_m =\frac{ (\ga_m +t_1+ t_0)(\ga_m+t_1-t_0)}{2\ga_m}-t_1.
\] 
\end{prop}
If we define $b_0=t_1$, which is consistent with the above defined expression for $b_m$, then the expression for $T_1p_{m}$, $m$ even, is also valid for $m=0$ if we assume that $p_{-1}$ exists (e.g.~$p_{-1}=1$). 
\begin{proof}
From $p_0=1$ we obtain directly $T_1 p_0 = t_1 p_0$. 
Assume $m \geq 1$. From the action of $T_1$ on the monomial $x^m$ (see the proof of Proposition \ref{prop:trian}) it follows that 
\[
T_1 p_{2n-1} = -p_{2n} + \sum_{j=0}^{2n-1} a_{2n-1,j} p_j, \qquad T_1 p_{2n} = \sum_{j=0}^{2n} a_{2n,j} p_j, \quad n \in \N.
\]
Now using $YT_1+T_1Y= Y^2+t_1^2-t_0^2$ we find $(Y+\ga_{m})T_1p_{m} = (\ga_{m}^2+t_1^2-t_0^2)p_{m}$. Applying this to the expansions gives
\[
a_{m,m}= \frac{ \ga_{m}^2+t_1^2-t_0^2 }{ 2\ga_{m} }, \qquad a_{2n-1,2n-2}=0,\qquad a_{m,j}= 0, \ j \leq m-2.
\]
Note that we do not find the value for $a_{2n,2n-1}$ from this, since $\ga_{2n-1}+\ga_{2n}=0$.
In order to find $a_{2n,2n-1}$ we write out $T_1^2p_{2n}$:
\[
\begin{split}
T_1^2p_{2n} &= T_1(a_{2n,2n}p_{2n}+a_{2n,2n-1}p_{2n-1})  \\
&=(a_{2n,2n}^2 - a_{2n,2n-1})p_{2n} + a_{2n,2n-1}( a_{2n,2n}+a_{2n-1,2n-1} ) p_{2n-1} .
\end{split}
\]
Using $\ga_{2n}=-\ga_{2n-1}$ we see that $a_{2n,2n}+a_{2n-1,2n-1}=0$.
Now from $T_1^2=t_1^2$ it follows that $(a_{2n,2n}^2 - a_{2n,2n-1})=t_1^2$, which gives us the value of $a_{2n,2n-1}$.
\end{proof}
We denote by $\mathcal A_Y \subset \mathcal T$ and $\mathcal A_Y^W \subset T$ the commutative subalgebras consisting of elements $p(Y)$ for $p \in \mathcal A$, respectively $p \in \mathcal A^W$. Also, we set $\mathcal A(n) = \mathrm{span}\{p_{2n}, p_{2n-1}\}$, $n \in \N$, and $\mathcal A(0) = \mathrm{span}\{p_0\}$. 
\begin{thm} \label{thm:decomp}
{\rm (a)} The representation $\pi_{\mathbf t} : \C[\mathcal W] \rightarrow \mathrm{End}(\mathcal A)$ is faithful.

{\rm (b)} The decomposition $\mathcal A = \bigoplus_{n \in \Z_{\geq 0}} \mathcal A(n)$ is the multiplicity-free, irreducible decomposition of $\mathcal A$ as a $\C[\mathcal W]$-module. Moreover, the element $p(\xi)$ in the center $\mathcal A_\xi^W$ of $\C[\mathcal W]$ acts on $\mathcal A(n)$ as multiplication by $p(\ga_{2n})$. 
\end{thm}
\begin{proof}
{\rm (a)} By Proposition \ref{prop:iso} every element of $\mathcal T=\pi_{\mathbf t}(\C[\mathcal W])$ can be written as $p(Y)+T_1q(Y)$ for some $p,q \in \mathcal A$. Suppose that $p(Y)+T_1q(Y)=0$ in $\mathcal T$, then $\big(p(\ga_m) + q(\ga_m)T_1\big) p_m =0$ for all $m\in \Z_{\geq 0}$. Now it follows from Proposition \ref{prop:V1} that $q(\ga_m)=0$ for all $m$, hence $q=0$ in $\mathcal A$. This leaves us with $p(\ga_m)=0$ for all $m$, so also $p=0$ in $\mathcal A$.

{\rm (b)} The decomposition follows from Proposition \ref{prop:V1} and the fact that the algebra $\mathcal T$ is generated by $T_1$ and $Y$. Let $p(\xi) \in \mathcal A^W_\xi$, then since $\ga_{2n}=-\ga_{2n-1}$ we find that $\pi_{\mathbf t}(p(\xi))=p(Y)$ acts on $\mathcal A(n)=span\{p_{2n-1},p_{2n}\}$ as multiplication by $p(\ga_{2n-1})=p(\ga_{2n})$.
\end{proof}

\subsection{A degenerate double affine Hecke algebra}
Let $z \in \mathrm{End}(\mathcal A)$ be multiplication by $x$; $(zp)(x)=xp(x)$ for $p \in \mathcal A$. We denote by $\mathcal A_z$ the commutative algebra in $\mathrm{End}(\mathcal A)$ of multiplication operators $p(z)$ for $p\in \mathcal A$. Let us now define a degenerate double affine Hecke algebra as follows. 
\begin{Def}
The algebra $\mathcal H=\mathcal H(\mathbf t)$ is the algebra in $\mathrm{End}(\mathcal A)$ generated by $\mathcal T$ and by $\mathcal A_z$.
\end{Def}
From writing the explicit expressions $T_i$, $i=0,1$, as
\[
T_i = t_i s_i + \frac{ t_i^2 - u_i^2 + \big(a_i^\vee(\cdot)\big)^2 }{ a_i(\cdot) }(s_i-1)
\]
it follows that we have the following relations in $\mathcal H$,
\begin{equation} \label{eq:relations}
T_i p(z) - (s_ip)(z) T_i = \frac{t_i^2-u_i^2+\big(a_i^\vee(z)\big)^2}{a_i(z)}\big( (s_i p)(z) - p(z) \big).
\end{equation}
Using the linear bases for $\C[\mathcal W]$ from Proposition \ref{prop:iso} we obtain linear bases for $\mathcal H$.
\begin{prop}
The sets $\{Y^kz^l ,Y^m T_1z^n \ |  \ k,l,m,n \in \Z_{\geq 0}\}$ and $\{ z^k Y^l, z^m T_1 Y^n\ | \ k,l,m,n \in \Z_{\geq 0} \}$ are linear bases for $\mathcal H$. So
\[
\mathcal A_Y \otimes T \otimes \mathcal A_z \cong \mathcal H \cong \mathcal A_z \otimes T \otimes \mathcal A_Y
\]
as vector spaces.
\end{prop}
Now we can show that the relations \eqref{eq:relations} and the quadratic relations in $\mathcal T$ completely characterizes the algebra $\mathcal H$.
\begin{prop} \label{prop:Ch-alg1}
The algebra $\mathcal H(\mathbf t)$ is isomorphic to the unital, associative algebra $\mathcal V(\mathbf t)$ generated by $V_0, V_1, v$ with relations
\[
V_i^2 =t_i^2, \qquad i=0,1,
\]
and, for $p \in \mathcal A$,
\[
V_i p(v) - (s_ip)(v) V_i = \frac{t_i^2-u_i^2+\big(a_i^\vee(v)\big)^2}{a_i(v)}\big( (s_i p)(v) - p(v) \big),
\]
The isomorphism $\phi:\mathcal V(\mathbf t) \rightarrow \mathcal H(\mathbf t)$ is given on generators by the assignments 
\[
V_0 \mapsto T_0, \quad V_1 \mapsto T_1, \quad v \mapsto z.
\]
\end{prop}
\begin{proof}
Let $\phi$ be defined on the generators of $\mathcal V$ as in the proposition, then by $T_i^2=t_i^2$, $i=0,1$, and by \eqref{eq:relations} $\phi$ preserves the defining relations of $\mathcal V$. Moreover, $\phi$ maps generators of $\mathcal V$ to generators of $\mathcal H$, hence $\phi$ is surjective. For the injectivity it is enough to show that the set 
\[
B=\{(V_0+V_1)^k z^l, (V_0+V_1)^m V_1 z^n \ | \ k,l,m,n \in \Z_{\geq 0}\}
\]
is a linear basis for $\mathcal V$. Now observe that the relations $V_i^2=t_i^2$ imply that the subalgebra generated by $V_0$ and $V_1$ is isomorphic to $\C[\mathcal W]$, then from Proposition \ref{prop:iso} and the other relations in $\mathcal V$ it follows that the set  $B$ is indeed a basis for $\mathcal V$.
\end{proof}

For $p(z)=z \in \mathcal A_z$ the relations \eqref{eq:relations} become
\[
\begin{split}
T_0 z + z T_0-T_0 &= t_0^2 - u_0^2 + (\hf-z)^2,\\ 
T_1 z + z T_1 &= u_1^2 - t_1^2 - z^2,
\end{split}
\]
which we can also write as
\[
\begin{split}
(T_0+\hf-z)^2=u_0^2, \qquad
(T_1+z)^2 =u_1^2.
\end{split}
\]
Let us therefore denote
\[
T_{a_0^\vee} = -T_0-\hf+z, \qquad T_{a_1^\vee} = -T_1-z,
\]
then we have operators associated to any simple root in $\mathcal R$ such that 
\[
T_{\al}^2 = t_\al^2, \qquad \al \in \{a_0,a_1,a_0^\vee,a_1^\vee\}.
\]
We will also use the notation $T_{a_i^\vee}=U_i$, $i=0,1$. Using Proposition \ref{prop:Ch-alg1} we now obtain the following characterization of the algebra $\mathcal H$.
\begin{prop} \label{prop:Ch-alg2}
The algebra $\mathcal H(\mathbf t)$ is isomorphic to the unital, associative algebra $\mathcal V(\mathbf t)$ generated by $V_i, \widetilde V_i$, $i=0,1$, with relations
\[
\begin{split}
&V_i^2=t_i^2, \quad \widetilde V_i^2 = u_i^2, \\
&V_0+V_1+\widetilde V_0 + \widetilde V_1 = -\hf.
\end{split}
\]
The isomorphism is given explicitly on generators by $V_i \mapsto T_i$, $\widetilde V_i \mapsto U_i$, for $i=0,1$.
\end{prop}
\begin{proof}
The last relation in the proposition implies that $\mathcal V$ is also generated as an algebra by the elements $V_0$, $V_1$ and $v=-V_1-\widetilde V_1=\hf+V_0+\widetilde V_0$. Writing out the quadratic relations $\widetilde V_i^2=u_i^2$ leads to 
\[
\begin{split}
V_0 v + v V_0-V_0 &= t_0^2 - u_0^2 + (\hf-v)^2,\\ 
V_1 v + v V_1 &= u_1^2 - t_1^2 - v^2.
\end{split}
\]
Iterating these relations gives back the defining relations from Proposition \ref{prop:Ch-alg1}.
\end{proof}

Finally, to justify the name ``degenerate double affine Hecke algebra'' for our algebra $\mathcal H$ we show that the operators $T_0$ and $T_0$ can be obtained by taking an appropriate limit in the Noumi representation of the affine Hecke algebra of type $\tilde A_1$. 

Let $0 <q<1$ and let $\mathbf k:\mathcal R \rightarrow \C$ be a multiplicity function. The affine Hecke algebra of type $\tilde A_1$ is the unital complex algebra generated by $V_0$ and $V_1$ with relations
\[
(V_j-k_{a_j})(V_j+k_{a_j}^{-1})=0, \qquad j=0,1.
\]
The Noumi representation is given by
\[
\begin{split}
V_j = k_{a_j} + k_{a_j}^{-1} \frac{ (1-k_{a_j} k_{a_j^\vee} q^{a_j^\vee(x)})(1-k_{a_j} k_{a_j^\vee}^{-1} q^{a_j^\vee(x)})}{1-q^{a_j(x)}}(s_j-1)
\end{split}
\]
These operators act on the algebra of Laurent polynomials in $q^x$. 
Now we substitute
\[
(k_{a_0}, k_{a_0^\vee}, k_{a_1} , k_{a_1^\vee}) \mapsto (-iq^{t_0}, iq^{u_0}, -iq^{t_1}, iq^{u_1}),
\]
with $t_0,u_0, t_1, u_1 \in \C$, then
\[
\lim_{q \uparrow 1} \frac{ 1-iV_j }{1-q} = t_j + \frac{ (t_j+u_j+a_j^\vee(x))(t_j-u_j+a_j^\vee(x))} {a_j(x)}(s_j-1)=T_j, \qquad j=0,1.
\]
So the algebra $\mathcal T$ generated by $T_0$ and $T_1$ is indeed a degenerate affine Hecke algebra. Moreover, the Dunkl-Cherednik operator for the affine Hecke algebra is given by $\widetilde Y=V_1V_0$. In the limit we obtain
\[
\lim_{q \uparrow 1} \frac{ \widetilde Y+1}{1-q} = \lim_{q \uparrow 1} \Big(\frac{ 1-iV_0}{1-q}+ \frac{ 1-iV_1}{1-q} - \frac{ (1-iV_0)(1-iV_1)}{1-q} \Big)= T_0+T_1 = Y.
\]

The algebra $\mathcal H$ may be considered as a $q=1$ analogue of the double affine Hecke algebra of type $(C^\vee_1, C_1)$, which was introduced by Sahi \cite{Sa99} for general rank. The presentation of $\mathcal H$ in Proposition \ref{prop:Ch-alg2} corresponds to Stokman's characterization \cite[Theorem 2.22]{NS04} of Sahi's double affine Hecke algebra.

\section{Wilson polynomials}\label{sec3}
In this section we study the non-symmetric Wilson polynomials from Definition \ref{def:Wilson pol}, and their symmetrized versions. We show that the symmetric Wilson polynomials are (with a suitable normalization) exactly the Wilson polynomials as defined by Wilson in \cite{Wi80}.

\subsection{Orthogonality relations} \label{ssec:orth rel}
Let us introduce parameters $a,b,c,d$ related to the multiplicity function $\mathbf t=(t_0,u_0,t_1,u_1)$ by
\begin{equation}  \label{eq:abcd}
(a,b,c,d) = (t_1+u_1, t_1-u_1,t_0+u_0+\hf, t_0-u_0+\hf).
\end{equation}
Throughout the rest of the paper the parameters $a,b,c,d$ will be related to the multiplicity function $\mathbf t$ in this way. 
We define a weight function $\De$ by
\[
\De(x;\mathbf t) = \frac{ \Ga(a + x)\Ga(a+1-x) \Ga(b + x) \Ga(b+1-x) \Ga( c\pm x) \Ga(d\pm x)}{ \Ga(2x)\Ga(1-2x)}.
\]
Here, and elsewhere, we use the notation $f(\al\pm \be)=f(\al+\be)f(\al-\be)$. From here on we assume that the multiplicity function $\mathbf t$ is such that
\begin{itemize}
\item $a,b,c,d \not\in -\hf \Z_{\geq 0}$,
\item the pairwise sum of $a,b,c,d$ is not contained in $\Z$.
\end{itemize}
Let $\mathcal C=\mathcal C_{\mathbf t}$ be a contour in the complex plane that runs along the imaginary axis from $-i \infty$ to $i\infty$ and is indented  such that the sequences $a+n$, $b+n$, $c+n$, $d+n$, $n \in \Z_{\geq 0}$, are separated by $\mathcal C$ from the sequences $-(a+n)$, $-(b+n)$, $-(c+n)$, $-(d+n)$, $n \in \Z_{\geq 0}$. Moreover, we assume that, set theoretically, $\mathcal C=-\mathcal C$. With the above assumptions on $\mathbf t$ such a contour exists. To the weight function $\De$ we now associate a non-degenerate bilinear form $\langle \cdot, \cdot \rangle_{\mathbf t}$ on $\mathcal A$,
\[
\langle f,g \rangle_{\mathbf t} = \frac{ 1}{2\pi i} \int_{\mathcal C} f(x) g(x)\De(x) dx.
\]
By Cauchy's Theorem we may write the above integral as an integral over $i\R$ plus a finite sum of residues corresponding to poles of $\De$. In case $a,b,c,d$ have positive real part we may take $\mathcal C=i\R$.
\begin{prop}
The operators $Y$ and  $T_i$, $i=0,1$, are symmetric with respect to $\langle \cdot, \cdot \rangle_{\mathbf t}$. 
\end{prop}
\begin{proof}
Let $f,g \in \mathcal A$. We have
\[
(T_{0}f)(x) = t_0f(x) + \frac{ (c-x)(d-x) }{1-2x}\big(f(1-x)-f(x)\big),
\]
so
\[
\langle T_{0}f, g \rangle_{\mathbf t} = \frac{ 1}{2\pi i} \int_{\mathcal C} t_0 f(x) g(x) \De(x) dx + \frac{1}{2\pi i} \int_{\mathcal C} \left( f(1-x) - f(x) \right ) g(x) \hat \De(x)dx,
\]
where
\[
\hat \De(x) = \frac{ \Ga(a + x)\Ga(a+1-x) \Ga(b + x) \Ga(b+1-x) \Ga( c+x)\Ga(c+1-x) \Ga(d+x) \Ga(d+1-x)}{ \Ga(2x)\Ga(2-2x)}.
\]
We see that $\hat \De(x) = \hat \De(1-x)$. So we may write the second integral as
\[
\frac1{2\pi i}\int_{1+\mathcal C} f(x)g(1-x) \hat \De(x)dx - \frac1{2\pi i}\int_{\mathcal C} f(x)g(x) \hat \De(x)dx.
\]
Under the current assumptions on $\mathbf t$ the integrand of the first integral does not have poles inside the area between $\mathcal C$ and $1+\mathcal C$. Therefore we can shift the path of integration, and we obtain
\[
\begin{split}
\langle T_{0}f, g \rangle_{\mathbf t}& = \frac{ 1}{2\pi i} \int_{\mathcal C} f(x) \left(t_0 g(x) + \frac{ (c-x)(d-x) }{1-2x}(g(1-x) - g(x)) \right) \De(x)dx \\
& = \langle f, T_0g \rangle_{\mathbf t}. 
\end{split}
\]
This proves the proposition for $T_0$. The proof for $T_1$ is similar. Since $Y=T_0+T_1$, $Y$ is also symmetric with respect to $\langle \cdot, \cdot \rangle_{\mathbf t}$.. 
\end{proof} 
We can now prove the orthogonality relations for the non-symmetric Wilson polynomials.
\begin{thm} \label{thm:orth}
The set of non-symmetric Wilson polynomials $\{p_m \ | \ m \in \Z_{\geq 0} \}$ is an orthogonal basis for $\mathcal A$ with respect to $\langle \cdot, \cdot\rangle_{\mathbf t}$.
\end{thm}
\begin{proof}
Using the previous proposition we have
\[
\ga_{m_1} \langle p_{m_1}, p_{m_2} \rangle_{\mathbf t} = \langle Yp_{m_1}, p_{m_2} \rangle_{\mathbf t} = \langle p_{m_1}, Yp_{m_2} \rangle_{\mathbf t} = \ga_{m_2} \langle p_{m_1}, p_{m_2} \rangle_{\mathbf t}.
\]
Since $\ga_{m_1} \neq \ga_{m_2}$ if $m_1 \neq m_2$, we have
\[
\langle p_{m_1}, p_{m_2} \rangle_{\mathbf t} = 0,
\]
if $m_1 \neq m_2$. 
\end{proof}
We evaluate the diagonal terms $\langle p_m, p_m \rangle_{\mathbf t}$ later on in Theorem \ref{thm:FG}.
\begin{rem}
We could also define a bilinear form by
\[
\langle f, g \rangle_{\mathbf t}' = \frac{1}{2\pi i} \int_{\mathcal C} f(x) g(-x) \De(x) dx.
\]
This bilinear form is closer to the bilinear form used for the non-symmetric Askey-Wilson polynomials in \cite{NS04}. However, $T_0$ and $T_1$ are not symmetric with respect to this bilinear form. So this would lead to biorthogonality relations between the eigenfunctions of $Y$ and $Y^*$, the adjoint of $Y$ with respect to the above defined bilinear form. However, it is not hard to show that the eigenfunctions of $Y^*$ are precisely the non-symmetric Wilson polynomials $x \mapsto p_m(-x;\mathbf t)$, so that this biorthogonality relation is equivalent to the orthogonality relations in Theorem \ref{thm:orth}.
\end{rem}

\subsection{Symmetric and anti-symmetric Wilson polynomials}
Let us define in $\C[W]$
\[
c_{\pm} = \hf(1\pm  s_1).
\]
These elements are orthogonal primitive idempotents in $\C[W]$, i.e., $c_\pm^2=c_\pm$, $c_-+c_+=1$ and $c_\pm c_{\mp}=0$. For $p \in \mathcal A$ the polynomial $c_+p$ is an even polynomial, so $c_+$ is the projection of $\mathcal A$ onto $\mathcal A^W$, and $c_-p$ is an odd polynomial. The representation $\pi_{\mathbf t}$ of $\C[W]$ gives us corresponding elements in the algebra $T$, 
\[
C_\pm = \pi_{\mathbf t}(c_{\pm}) = \frac{ 1}{2t_1}(t_1 \pm T_1) \in  T.
\]
Since $C_++C_-=1$ in $T$, we have a corresponding decomposition of $\mathcal A$ in a symmetric and an anti-symmetric part; $\mathcal A = \mathcal A_+ \oplus \mathcal A_-$, where $\mathcal A_\pm = C_\pm \mathcal A$. So $\mathcal A_\pm$ consists of polynomials $p \in \mathcal A$ for which $(T_1 \mp t_1) p = 0$. From the explicit expression of $T_1$ we obtain $T_1p=t_1p$ if and only if $s_1p=p$, so we see that $\mathcal A_+ = \mathcal A^W$, the algebra of even polynomials. Moreover, $\mathcal A_-$ consists of the polynomials $p$ such that $c_1(-x)p(x)$ is odd. Indeed, the identity $(T_1+t_1)p =0$ is equivalent to
\[
0=2t_1p(x) + c_1(x)\big(p(-x)-p(x)\big) = \big(c_1(x)+c_1(-x)\big)p(x) + c_1(x)\big(p(-x)-p(x)\big),
\]
which gives $c_1(-x)p(x)=-c_1(x)p(-x)$. The symmetric and the anti-symmetric polynomials are related by the generalized Weyl denominator, which is the monic anti-symmetric polynomial of lowest degree.
\begin{prop} \label{prop:A+ A-}
Let the generalized Weyl denominator $\de$ be the polynomial given by 
\[
\de(x)=(t_1+u_1+x)(t_1-u_1+x),
\]
then $\de(z) \mathcal A_+ = \mathcal A_-$. 
\end{prop}
\begin{proof}
After a straightforward calculation it follows from Proposition \ref{prop:Ch-alg1} that 
\[
(T_1 + t_1)\de(z) = \de(-z)(T_1-t_1),
\]
hence $\de(z) \mathcal A_+ \subset \mathcal A_-$. In the same way it follows that $\de(z)^{-1} \mathcal A_-$ is $W$-invariant. So we only need to show that if $p \in \mathcal A_-$, then $\de(x)^{-1}p(x)$ is a polynomial. This follows directly from writing the polynomial $p  \in \mathcal A_-$ as
\[
p(x) =  -\frac{ c_1(x) }{c_1(-x)}p(-x) = \frac{ \de(x) }{\de(-x)}p(-x),
\]
where we used $\de(x)=2xc_1(x)$. Hence $\de(z)^{-1}\mathcal A_- \subset \mathcal A_+$.
\end{proof}
We can decompose the irreducible $\mathcal T$-modules in a symmetric and anti-symmetric part; $\mathcal A(n) = \mathcal A_+(n) \oplus \mathcal A_-(n)$, $n \in \Z_{\geq 0}$. So the non-symmetric Wilson polynomials $p_m$ can be written as the sum of a symmetric and an anti-symmetric polynomial. For a factor $F^\pm_n$ independent of $x$, we must have $C_\pm p_{2n-1} = F_n^{\pm} C_\pm p_{2n}$, since both $C_\pm p_{2n}$ and $C_\pm p_{2n-1}$ are in the (anti-)symmetric part of $\mathcal A(n)$. It will be useful to work with monic polynomials. 
\begin{Def}
{\rm (a)} The symmetric Wilson polynomial $P_{2n}^+$ is the unique monic polynomial in $\mathcal A_+ (n)$. 

{\rm (b)} The anti-symmetric Wilson polynomial $P_{2n}^-$ is the unique monic polynomial in $\mathcal A_- (n)$.
\end{Def}
We remark that both $P_{2n}^+$ and $P_{2n}^-$ are of degree $2n$. There are no (anti-)symmetric Wilson polynomials of odd degree. Note that the sets $\{P_{2n}^+\ | \ n \in \Z_{\geq 0} \}$ and $\{P_{2n}^- \ | \ n \in \N\}$ are a linear basis for $\mathcal A_+$, respectively $\mathcal A_-$. We can express the (anti-)symmetric Wilson polynomials in terms of non-symmetric ones and vice versa.
\begin{lem} \label{lem:sym pol}
We have 
\[
P_{2n}^+ = p_{2n} + (b_{2n}-t_1)p_{2n-1} , \qquad
P_{2n}^- = p_{2n} + (b_{2n}+t_1)p_{2n-1},
\]
or equivalently,
\[
p_{2n}=\frac{1}{2t_1}\Big((b_{2n}+t_1)P_{2n}^+ -(b_{2n}-t_1)P_{2n}^-\Big), \qquad p_{2n-1}= \frac{1}{2t_1}\Big(P_{2n}^- -P_{2n}^+\Big),
\]
with $b_{2n}$ as in Proposition \ref{prop:V1}.
\end{lem}
\begin{proof}
This follows from writing out $(T_1 \pm t_1)p_m$ using Proposition \ref{prop:V1}.
\end{proof}

The (anti-)symmetric Wilson polynomials also satisfy orthogonality relations with respect to the bilinear form $\langle\cdot,\cdot\rangle_{\mathbf t}$.
\begin{lem} \label{lem:orth2}
We have the following orthogonality relations:

{\rm (a)} For $n \in \Z_{\geq 0}$, $m \in \N$, $\langle P^+_{2n},P^-_{2m} \rangle_{\mathbf t}=0$, 

{\rm (b)} For $n,m \in \Z_{\geq 0}$, $n \neq m$, $\langle P^+_{2n}, P_{2m}^+ \rangle_{\mathbf t}=0$,

{\rm (c)} For $n,m \in \N$, $n\neq m$, $\langle P_{2n}^-, P_{2m}^- \rangle_{\mathbf t}=0$.
\end{lem}
\begin{proof}
{\rm (a)} From $C_\pm C_\mp=0$ we obtain
\[
\langle C_+p_n, C_-p_m \rangle_{\mathbf t} = \langle p_n, C_+C_-p_m \rangle_{\mathbf t} =0.
\]

{\rm (b)} Since $P_{2n}^+ \in \mathcal A^W$ and $Y^2 \in \mathcal A_Y^W$, it follows from Theorem \ref{thm:decomp}{\rm (b)} that $Y^2P_{2n}^+=\ga_{2n}^2 P_{2n}^+$. Because $Y$ is symmetric with respect to $\langle \cdot,\cdot\rangle_{\mathbf t}$, $Y^2$ is also symmetric. Now, since the eigenvalues $\ga_{2n}^2$ are pairwise different, the orthogonality relations follow.

{\rm (c)} Theorem \ref{thm:orth} gives for $n \neq m$ 
\[
0=\langle p_n, p_m \rangle_{\mathbf t} = \langle C_+p_n, C_+p_m \rangle_{\mathbf t}+\langle C_+p_n, C_-p_m \rangle_{\mathbf t}+\langle C_-p_n, C_+p_m \rangle_{\mathbf t}+\langle C_-p_n, C_-p_m \rangle_{\mathbf t}.
\]
Then by {\rm (a)} and {\rm (b)} we obtain $\langle C_-p_n, C_-p_m \rangle_{\mathbf t}=0$.
\end{proof}

Let us define the weight function $\De^+$ by
\[
\De^+(x;\mathbf t) = \frac{ \Ga(a \pm x) \Ga(b \pm x) \Ga(c \pm x) \Ga(d \pm x) }{ \Ga( \pm 2x )},
\]
and let $\langle \cdot, \cdot \rangle_{\mathbf t}^+$ be the corresponding bilinear form on $\mathcal A_+$,
\[
\langle f, g \rangle_{\mathbf t}^+ = \frac{1}{2\pi i} \int_{\mathcal C} f(x) g(x) \De^+(x) dx.
\]
Observe that $\De^+(-x)= \De^+(x)$.
\begin{lem} \label{lem:blforms}
For $f,g \in \mathcal A_+$
\[
\langle f, g \rangle_{\mathbf t} = \hf(a+b) \langle f, g \rangle_{\mathbf t}^+.
\]
\end{lem}
\begin{proof}
From the explicit expression for $\De$ and $\De^+$ it follows that
\[
\De(x) = c_1(-x)\De^+(x).
\]
Let $f,g \in \mathcal A_+$, then 
\[
2\langle f,g \rangle_{\mathbf t} =\frac{1}{2\pi i}\int_{\mathcal C}f(x) g(x) \big( c_1(x) + c_1(-x) \big) \De^+(x) dx.
\]
Using $c_1(x) + c_1(-x) = 2t_1=a+b$, we obtain
\[
\langle f,g \rangle_{\mathbf t} = \hf(a+b)\langle f,g \rangle_{\mathbf t}^+. \qedhere
\]
\end{proof}
Combining Lemma \ref{lem:blforms} with Lemma \ref{lem:orth2}{\rm (b)} then leads to the orthogonality relations for the symmetric Wilson polynomials with respect to $\langle\cdot, \cdot \rangle_{\mathbf t}^+$.
\begin{thm} \label{thm:orth+}
The set $\{P_{2n}^+\ | \ n \in \Z_{\geq 0}\}$ is an orthogonal basis for $\mathcal A_+$ with respect to $\langle \cdot,\cdot \rangle_{\mathbf t}^+$.
\end{thm}
Now we can prove the generalized Weyl character formula, which says that anti-symmetric Wilson polynomials can be expressed in terms of symmetric Wilson polynomials with a shift in the parameters using the generalized Weyl denominator $\de$ defined in Proposition \ref{prop:A+ A-}.
\begin{thm} \label{thm:shift}
For $n \in \N$
\[
P^-_{2n}(x;\mathbf t) = \de(x) P^+_{2n-2}(x;t_0, u_0, t_1+1, u_1).
\]
\end{thm}
\begin{proof}
The symmetric Wilson polynomial $P^+_{2n}(\cdot;t_0,u_0,t_1+1,u_1)$ is the unique even monic polynomial of degree $2n$ that is orthogonal to all even polynomials of degree $\leq 2n-2$ with respect to $\langle \cdot,\cdot\rangle^+_{(t_0,u_0,t_1+1,u_1)}$. Let $p(x) = \de(x)^{-1} P_{2n}^-(x;\mathbf t)$, then $p$ is a monic polynomial and by Proposition \ref{prop:A+ A-} we have $p \in \mathcal A_+$. So, to prove the desired identity, it is enough to show that $\langle p, q_{2k} \rangle_{(t_0,u_0,t_1+1,u_1)}^+=0$, for any even polynomial $q_{2k}$ of degree $2k$ for $k =0,\ldots,n-2$.

By the explicit expressions for $\De$ and $\De^+$ we have
\[
\de(x) \De(x;\mathbf t) = -\frac{1}{2x} \De^+(x;t_0,u_0,t_1+1,u_1).
\]
Now we write out the bilinear form $\langle \cdot,\cdot \rangle_{\mathbf t}$ as an integral and we symmetrize the integrand, then we obtain
\[
\begin{split}
\langle P_{2n}^-, \de(z)q_{2k} \rangle_{\mathbf t} &= \frac{1}{4\pi i} \int_{\mathcal C} p(x)q_{2k}(x) \Big(\frac{ \de(-x) - \de(x) }{2x} \Big) \De^+(x;t_0,u_0,t_1+1,u_1) dx\\
&= -\hf(a+b) \langle p, q_{2k} \rangle^+_{(t_0,u_0,t_1+1,u_1)}.
\end{split}
\]
Since $\de(z)q_{2k} \in \mathrm{span}\{ P^-_{2m}\ | \ m=1,\ldots, k+1\}\subset \mathcal A_-$ it follows from Lemma \ref{lem:orth2}{\rm (c)} that $\langle P_{2n}^-, \de(z)q_{2k} \rangle_{\mathbf t}=0$, hence $\langle p, q_{2k} \rangle^+_{(t_0,u_0,t_1+1,u_1)}=0$ as desired.
\end{proof}

\subsection{Duality}
In this subsection we prove the duality property for the Wilson polynomials using ideas from \cite{Sa99} and \cite{NS04}. We define an involution $\si$ acting on multiplicity functions by interchanging the values on the $a_0$-orbit and the $a_1^\vee$-orbit. So, given a multiplicity function $\mathbf t=(t_0,u_0,t_1,u_1)$, the multiplicity function $\mathbf t^\si:\mathcal R \rightarrow \C$ is given by
\[
\mathbf t^\si = (t_{a_0}^\si, t_{a_0^\vee}^\si, t_{a_1}^\si, t_{a_1^\vee}^\si) = (u_1, u_0, t_1, t_0).
\]
We call $\mathbf t^\si$ the dual of $\mathbf t$. If an object depends on the multiplicity function $\mathbf t$, we will use a super- or subscript $\si$ to denote the same object depending on $\mathbf t^\si$. For instance, for the difference-reflection operators we write $T_i^\si = T_i(\mathbf{t}^\si)$, and $\mathcal H_\si=\mathcal H(\mathbf t^\si)$ is the algebra generated by $T_0^\si$, $T_1^\si$ and $z^\si=z$. 
\begin{prop} \label{prop:dual iso}
The assignments 
\[
T_{0} \mapsto -(T_1^\si+ z), \qquad T_1 \mapsto T_1^\si, \qquad z \mapsto -(T_0^\si + T_1^\si),
\]
extend uniquely to an algebra isomorphism $\si : \mathcal H \rightarrow \mathcal H_\si$ with inverse $\si^{-1}=\si_\si$, respectively anti-isomorphism $\psi : \mathcal H \rightarrow \mathcal H_\si$ with inverse $\psi^{-1}=\psi_\si$.
\end{prop}
We call $\si$ and $\psi$, the duality isomorphism, respectively duality anti-isomorphism. 
\begin{proof}
It is straightforward to verify that the following relations are satisfied
\[
\si(T_0)^2=t_0^2, \quad \si(T_1)^2=t_1^2, \quad \si(z-\hf-T_0)^2 = u_0^2, \quad \si(-T_1-z)^2=u_1^2,
\]
and similarly for $\psi$. We see that both $\psi$ and $\si$ satisfy
\begin{gather*}
T_{0}\mapsto U_1^\si, \qquad T_{1} \mapsto T_{1}^\si, \qquad U_0 \mapsto  U_0^\si, \qquad U_1 \mapsto T_{0}^\si, \qquad
z \mapsto -  Y^\si, \qquad Y \mapsto -z.
\end{gather*}
We then see that $\si^{-1}= \si_\si$ and $\psi^{-1} = \psi_\si$.
\end{proof}

Next we introduce, following Sahi \cite{Sa99}, elements in $\mathcal H$ called intertwiners which can be used to construct raising and lowering operators for the non-symmetric Wilson polynomials. These intertwiners $S_0, S_1 \in \mathcal H$ are defined by
\begin{equation} \label{def:intertwiners}
S_0 = U_1 Y - Y U_1, \qquad S_1=T_1 Y - Y T_1.
\end{equation}
We have the following useful property.
\begin{lem} \label{lem:intertwiners}
For $n \in \Z_{\geq 0}$,
\[
S_0 p_{2n} = (\ga_{2n+1}-\ga_{2n}) p_{2n+1}, \qquad
S_1p_{2n+1} = 2\ga_{2n+2} p_{2n+2}.
\]
\end{lem}
\begin{proof}
By Proposition \ref{prop:V1} we have
\[
S_1p_{2n+1} = (\ga_{2n+1}-Y) T_1 p_{2n+1} = (\ga_{2n+1}-Y)(-p_{2n+2} + b_{2n+1}\,p_{2n+1})
\]
Then the formula for $S_1$ follows from $\ga_{2n+2} = -\ga_{2n+1}$.

Next we check the action of $S_0$. Let $\widetilde S_0 = zT_0-T_0z \in \mathcal H$. For $p \in \mathcal A$, we obtain from Proposition \ref{prop:Ch-alg1} $p(z)\widetilde S_0 = \widetilde S_0 (s_0p)(z)$ in $\mathcal H$. Applying the duality isomorphism $\si$ and replacing the parameters by their duals, gives us $p(-Y) S_0 = S_0 p(1+Y)$ in $\mathcal H$. So $S_0p_{2n}$ is an eigenfunction of $Y$ for eigenvalue $-(1+\ga_{2n})=\ga_{2n+1}$, therefore $S_0p_{2n}=k_{n}p_{2n+1}$ for some constant $k_n$. Using Proposition \ref{prop:V1} again, the constant $k_n$ can be determined by finding the leading coefficient in
\[
S_0p_{2n} = (Y-\ga_{2n})(T_1+z)p_{2n} =(Y-\ga_{2n})((b_{2n}+z) p_{2n} + (b_{2n}^2-t_1^2)p_{2n-1}).
\]
From Proposition \ref{prop:trian} it now follows that $k_n = \ga_{2n+1}-\ga_{2n}$.
\end{proof}
Observe that in the same way as in Lemma \ref{lem:intertwiners} it can be proved that $S_0p_{2n+1}=k^0_n p_{2n}$ and $S_1p_{2n}=k^1_np_{2n-1}$ for some constants $k_n^i$ which can be determined explicitly. We do not need these formulas here.

As a first application of Lemma \ref{lem:intertwiners} we deduce a Rodriguez-type formula for the non-symmetric Wilson polynomials, which says 
that the non-symmetric Wilson polynomials can be generated from $1 \in \mathcal A$, the polynomial identically equal to $1$, using the intertwiners.
We use the standard notation for shifted factorials,
\[
(\al)_0=1, \qquad (\al)_n=\al(\al+1)\ldots(\al+n-1), \ n \in \N.
\]
\begin{prop} \label{prop:Rodriguez}
For $n \in \Z_{\geq 0}$,
\[
(S_1S_0)^n 1 = (-1)^n (2t_0+2t_1+1)_{2n}\, p_{2n}, \qquad S_0(S_1S_0)^n 1 = (-1)^{n+1} (2t_0+2t_1+1)_{2n+1}\,p_{2n+1}.
\]
\end{prop}

\begin{proof}
This follows from Lemma \ref{lem:intertwiners} using induction on the degree.
\end{proof}
As a second application of Lemma \ref{lem:intertwiners} we evaluate $p_m(-x_0)$ explicitly, where $x_0=t_1+u_1$. For this we introduce the evaluation mapping $\mathrm{Ev}:\mathcal H \rightarrow \C$ by
\[
\mathrm{Ev}(X) = \big(X(1)\big)(-x_0), \qquad X \in \mathcal H,
\]
where $1 \in \mathcal A$. From the explicit expression $T_1=t_1+c_1(\cdot)(s_1-1)$ and from $c_1(-x_0)=0$ it follows directly that
\begin{equation} \label{eq:Ev T1}
\mathrm{Ev}(T_1X)=t_1\mathrm{Ev}(X), \qquad X \in \mathcal H.
\end{equation}
\begin{prop} \label{prop:evaluation}
For $n \in \N$,
\[
\begin{split}
p_{2n}(-x_0) &= \frac{ (a+b+1)_n (a+c)_n (a+d)_n}{(n+a+b+c+d)_n},\\
p_{2n-1}(-x_0)&= - \frac{ (a+b+1)_{n-1} (a+c)_n (a+d)_n }{ (n+a+b+c+d-1)_n }, \\
P_{2n}^+(x_0) &= \frac{ (a+b)_n (a+c)_n (a+d)_n}{(n+a+b+c+d-1)_n}
\end{split}
\]
\end{prop}
\begin{proof}
Let $n \in \Z_{\geq 0}$. Using the relation $YT_1+T_1Y = Y^2+t_1^2-t_0^2$ in $\mathcal H$, we may write
\[
S_1 = 2T_1 Y - Y^2 -t_1^2+t_0^2.
\]
Then using \eqref{eq:Ev T1} we have
\[
\begin{split}
\mathrm{Ev}\big(S_1p_{2n+1}(z)\big) &= \big(2t_1 \ga_{2n+1}-\ga_{2n+1}^2-t_1^2 +t_0^2\big) p_{2n+1}(-x_0) \\
&= (t_0+t_1-\ga_{2n+1})(t_0-t_1+\ga_{2n+1})p_{2n+1}(-x_0).
\end{split}
\]
Applying the duality isomorphism $\si$ to the relation $T_0z+zT_0-T_0=t_0^2-u_0^2+(\hf-z)^2$ and replacing the parameters by their duals, we find in $\mathcal H$
\[
U_1Y+YU_1+U_1=u_0^2-u_1^2-(\hf+Y)^2.
\]
This gives us
\[
S_0=-2(T_1+z)Y-(T_1+z)+(\hf+Y)^2+u_1^2-u_0^2, 
\]
and using \eqref{eq:Ev T1} we then obtain
\[
\begin{split}
\mathrm{Ev}\big(S_0p_{2n}(z)\big) &= \big(-2(t_1-x_0)\ga_{2n}-(t_1-x_0)+(\hf+\ga_{2n})^2+u_1^2-u_0^2\big) p_{2n}(-x_0)\\
&=(u_1+u_0+\hf+\ga_{2n})(u_1-u_0+\hf+\ga_{2n})p_{2n}(-x_0).
\end{split}
\]
On the other hand, we find from Lemma \ref{lem:intertwiners}  
\[
\mathrm{Ev}\big(S_0p_{2n}(z)\big) = (\ga_{2n+1} - \ga_{2n})p_{2n+1}(-x_0), \qquad
\mathrm{Ev}\big(S_1p_{2n+1}(z)\big) = 2\ga_{2n+2} p_{2n+2}(-x_0).
\]
Combining now gives us the recurrence relations
\[
\begin{split}
p_{2n+1}(-x_0) &= -\frac{(t_0+t_1+u_0+u_1+\hf+n)(t_0+t_1-u_0+u_1+\hf+n)}{2t_0+2t_1+2n+1} p_{2n}(-x_0),\\
p_{2n+2}(-x_0) &= -\frac{(2t_0+2t_1+1+n)(2t_1+1+n)}{2t_0+2t_1+2n+2} p_{2n+1}(-x_0).
\end{split}
\]
The evaluation formula for $p_m(-x_0)$ follows by induction on $m$, starting with $p_0(-x_0)=1$. The expression for $P_{2n}^+(x_0)=P_{2n}^+(-x_0)$ follows from Lemma \ref{lem:sym pol}.
\end{proof}

Similar to $\mathrm{Ev}$ we define the dual evaluation mapping $\widetilde{\mathrm{Ev}}: \mathcal H_\si \rightarrow \C$ by
\[
\widetilde{\mathrm{Ev}}(\widetilde X) = \big(\widetilde X(1)\big)(-\ga_0), \qquad \widetilde X \in \mathcal H_\si.
\]
With the evaluation mappings and the duality anti-isomorphism $\psi$ we construct two pairings
$B: \mathcal H \times \mathcal H_\si \rightarrow \C$ and $\widetilde B: \mathcal H_\si \times \mathcal H \rightarrow \C$ as follows
\[
B(X,\widetilde X)= \mathrm {Ev}(\psi_\si(\widetilde X) X), \qquad \widetilde B(\widetilde X, X) = \widetilde { \mathrm {Ev}}(\psi(X) \widetilde X), \qquad X \in \mathcal H,\ \widetilde X \in \mathcal H_\si.
\]
These pairings have the following properties.
\begin{lem} \label{lem:pairing}
Let $X, X_1, X_2 \in \mathcal H$ and $\widetilde X, \widetilde X_1, \widetilde X_2 \in \mathcal H_\si$, and let $p \in \mathcal A$. Then

{\rm (a)} $B(X,\widetilde X) = \widetilde B(\widetilde X, X)$,

{\rm (b)} $B(X_1 X_2, \widetilde X) = B(X_2, \psi(X_1) \widetilde X)$, and $B(X, \widetilde X_1 \widetilde X_2) = B(\psi_\si(\widetilde X_1) X, \widetilde X_2)$,

{\rm (c)} $B\big( (Xp)(z), \widetilde X \big) = B\big( X \, p(z), \widetilde X)$, and $B\big( X, (\widetilde X p)(z)\big) = B\big( X , \widetilde X p(z) \big)$.
\end{lem}
\begin{proof}
{\rm (a)} Let $X = f(z)T_1g(Y) \in \mathcal H$ for $f,g \in \mathcal A$, then using $Y1= \ga_0$ and $T_1 1=t_1$ we have
\[
\begin{split}
\widetilde{\mathrm{Ev}}(\psi(X))\, 
&= \big(g(-z) T_1^\si f(-Y^\si) (1)\big)(-\ga_0) \\
&= t_1f(-x_0)g(\ga_0) \\
&= \big(f(z)T_1g(Y)(1)\big)(-x_0)\\
& = \mathrm{Ev}(X),
\end{split}
\] 
and similarly for $X=f(z)g(Y)$. So we have $\widetilde {\mathrm{Ev}}(\psi(X)) = \mathrm{Ev}(X)$ for all $X \in \mathcal H$. This gives us
\[
\widetilde B(\widetilde X, X) = \widetilde { \mathrm {Ev}}\big(\psi(X)\, \psi(  \psi_\si(\widetilde X)) \big) = \mathrm{Ev}\big(\psi_\si(\widetilde X) X \big) = B(X, \widetilde X).
\]

{\rm (b)} Since $\psi$ is an anti-isomorphism with inverse $\psi_\si$, we have
\[
B(X_1 X_2, \widetilde X) = \mathrm{Ev}\big(\psi_\si(\widetilde X)\, X_1 X_2\big) = \mathrm{Ev}\big(\psi_\si(\psi(X_1)\,\widetilde X)\, X_2\big) = B(X_2, \psi(X_1) \widetilde X).
\]

{\rm (c)} This is an immediate consequence of $\big((X p)(z)\big)(1) = (Xp) =  X\big(p(z)(1)\big)$ in $\mathcal A$.
\end{proof}
We renormalize the non-symmetric Wilson polynomials as follows:
\[
E(x,\ga_m;\mathbf t) = \frac{ p_m(x;\mathbf t)}{p_m(-x_0; \mathbf t)}.
\]
In particular we have $\mathrm{Ev}\big( E(z, \ga_m)\big)=1$ for all $m \in \Z_{\geq 0}$. Furthermore, we denote the eigenvalues of $Y^\si \in \mathcal H_\si$ by $x_m$, i.e., for $m \in \Z_{\geq 0}$,
\[
x_m=\ga_m^\si = 
\begin{cases}
t_1+u_1+n, & m=2n,\\
-(t_1+u_1+n), & m=2n-1.
\end{cases}
\]
We are now ready to prove the duality property for the non-symmetric Wilson polynomials.
\begin{thm} \label{thm:duality}
For $m,n \in \Z_{\geq 0}$ and $p \in \mathcal A$, we have
\[
p(-\ga_m) = \widetilde B\big(p(z), E(z,\ga_m) \big), \qquad  p(-x_n)= B\big( p(z), E_\si(z,x_n)\big) .
\]
Consequently, the Wilson polynomials satisfy the duality property
\[
E(-x_n, \ga_m;\mathbf t) = E(-\ga_m,x_n;\mathbf t^\si).
\]
\end{thm}
\begin{proof}
Using $\mathrm{Ev}\big( E(z, \ga_m)\big)=1$ for all $m \in \Z_{\geq 0}$ and the previous lemma, we have for $f \in \mathcal A$, 
\[
\begin{split}
\widetilde B\big(f(z), E(z, \ga_m) \big) \,&= \widetilde B\big( 1, f(-Y)\, E(z, \ga_m) \big)\\
& = \widetilde B\Big( 1, \big(f(-Y) E(\cdot, \ga_m)\big)(z) \Big)  \\
&=f(-\ga_m)\widetilde B\big( 1, E(z, \ga_m) \big)\\
&=f(-\ga_m).
\end{split}
\]
Similarly, for $g \in \mathcal A$, we find
\[
B\big( g(z), E_\si(z,x_n)\big) = g(-x_n).
\]
Now take $f= E_\si(\cdot;x_n)$ and $g=E(\cdot,\ga_m)$, and use Lemma \ref{lem:pairing}{\rm (a)} to obtain the second statement of the theorem.
\end{proof}
Theorem \ref{thm:duality} can be used to write the actions of $T_1$ and $U_1$ on the Wilson polynomials as difference-reflection operators in the dual variable. This will be useful in the next section.
\begin{lem} \label{lem:h1}
For $m \in \Z_{\geq 0}$
\[
\begin{split}
\big(T_1 E(\cdot;\ga_m)\big)(x) &= t_1 E(x;\ga_m) + c_{1}^\si(-\ga_m)\big(E(x;-\ga_m) - E(x;\ga_m)\big),\\
\big(U_1 E(\cdot;\ga_m)\big)(x) &= u_1 E(x;\ga_m) + c_{0}^\si(-\ga_m)\big(E(x;-1-\ga_m) - E(x;\ga_m)\big),
\end{split}
\]
where we use the convention $E(x,-\ga_0)=1$.
\end{lem}
It is useful to observe that $c_1^\si(-\ga_0)=0$.
\begin{proof}
For $U_1$ we recall that by Theorem \ref{thm:duality}
\[
\big(U_1 E(\cdot,\ga_m)\big)(-x_n) = B\big(E(z,\ga_m), T_{0}^\si E_\si(z, x_n) \big).
\]
Now we use the relations between $T_0$ and $p(z)$ in $\mathcal H$ as well as the identity $B(X,\widetilde X T_0^\si) = u_1 B(X, \widetilde X)$ for $X \in \mathcal H$, $\widetilde X \in \mathcal H_\si$, then we obtain for $m,n \in \Z_{\geq 0}$,
\[
\begin{split}
\big(U_1 & E(\cdot,\ga_m)\big)(-x_n) = u_1 B\big(E(z,\ga_m),  E_\si(1-z, x_n) \big) \\
&+ \frac{  u_1^2 - u_0^2 + (\hf+\ga_m)^2/4 }{1+2\ga_m} \Big( B\big(E(z,\ga_m),  E_\si(1-z, x_n) \big) - B\big(E(z,\ga_m),E_\si(z, x_n) \big) \Big)\\
& = u_1E(-x_n,\ga_m) + c_{0}(-\ga_m;\mathbf{t}^\si) \big(E(-x_n,-1-\ga_m)- E(-x_n,\ga_m) \big).
\end{split}
\]
Here we used
\[
B\big(E(z,\ga_m), E_\si(1-z, x_n) \big) = \widetilde B\big(E_\si(1- z,x_n), E(z, \ga_m) \big) = E_\si(1+\ga_m,x_n) = E(-x_n,-1-\ga_m).
\]
So the desired identity holds for all $n \in \Z_{\geq 0}$, hence it holds in $\mathcal A$. The expression for $T_1$ can be obtained in the same way.
\end{proof}

Let us also briefly consider the symmetric Wilson polynomials.
We denote $C_\pm E(\cdot,\ga_{m}) = E^\pm (\cdot,\ga_{m})$. 
For $m=2n$ or $m=2n-1$ we have $C_+E(\cdot,\ga_{m})=k P_{2n}^+$ for some constant $k$. From $c_1(-x_0)=0$ and $E(-x_0,\ga_{m})=1$ it follows that $(T_1 E(\cdot,\ga_{m}))(-x_0)=t_1$. Therefore $(C_+E(\cdot,\ga_{m}))(-x_0)=1$, which gives 
\[
E^+(x,\ga_{2n};\mathbf t) = E^+(x,\ga_{2n-1};\mathbf t)=\frac{ P^+_{2n}(x;\mathbf t) }{ P^+_{2n}(x_0;\mathbf t)}, \qquad n \in \N.
\]
In the same way as in Theorem \ref{thm:duality} the duality property for the renormalized symmetric Wilson polynomials is obtained;
\[
E^+(x_{2n}, \ga_{2m};\mathbf t) = E^+(\ga_{2m},x_{2n};\mathbf t^\si), \qquad m,n \in \Z_{\geq 0}.
\]

\subsection{Explicit expressions for the Wilson polynomials}
From the definition of $T_0$ and $T_1$, and from the description of $\mathcal W$ as the semi-direct product $W \ltimes \tau(\Z)$, it follows that we can write any $X \in \mathcal T$ as
\[
X = \sum_{\substack{\la \in \Z \\ w \in W}} c_{\la,w} \tau(\la) w
\]
for some coefficients $c_{\la,w} \in \C(x)$ of which only finitely many are non-zero. Since $W=\{1,s_1\}$, we may write any $X \in \mathcal T$ as $X=X_0 + X_1 s_1$, where $X_0,X_1 \in \bigoplus_{\la \in \Z} \C(x)\tau(\la)$. So $X_0$ and $X_1$ are difference operators with rational coefficients. We define $X_{sym} = X_0 + X_1$. Since $\mathcal A_+$ is $W$-invariant, it is clear the actions of $X$ and $X_{sym}$ on $\mathcal A_+$ coincide.

We use the difference operator $(Y^2)_{sym}$ to obtain a difference equation for the symmetric Wilson polynomials.
\begin{prop} \label{prop:diff eq}
{\rm (a)} The symmetric Wilson polynomials satisfy the equation
\[
L E^+(\cdot,\ga_{2n}) = n(n+a+b+c+d-1) E^+(\cdot,\ga_{2n}),
\]
where $L$ is the second-order difference operator
\begin{gather*}
L = (Y^2)_{sym}-(t_0+t_1)^2 =A(x)\Big(\tau(1)-1\Big) + A(-x)\Big(\tau(-1)-1\Big),\\
A(x) = \frac{ (a+x)(b+x)(c+x)(d+x) }{ 2x(2x+1) }.
\end{gather*}

{\rm (b)} The symmetric Wilson polynomials satisfy the recurrence relation
\[
(x^2-a^2)E^+(x,\ga_{2n})= B_n \Big(E^+(x,\ga_{2n+2})-E^+(x,\ga_{2n})\Big) +C_n \Big(E^+(x,\ga_{2n-2})-E^+(x,\ga_{2n})\Big),
\]
where
\[
\begin{split}
B_n &= \frac{ (n+a+b+c+d-1)(n+a+b)(n+a+c)(n+a+d) }{  (2n+a+b+c+d-1)(2n+a+b+c+d) },\\
C_n &= \frac{ n(n+b+c-1)(n+b+d-1)(n+c+d-1) }{ (2n+a+b+c+d-2)(2n+a+b+c+d-1) }.
\end{split}
\]
\end{prop}
Here we use the convention $E^+(x,\ga_{-2})=0$.
\begin{proof}
{\rm (a)} Since $E^+(\cdot,\ga_{2n}) \in \mathcal A(n)$ it follows from Theorem \ref{thm:decomp} that $Y^2 E^+(\cdot,\ga_{2n}) = \ga_{2n}^2 E^+(\cdot,\ga_{2n})$. Then, since $E^+(\cdot,\ga_{2n}) \in \mathcal A_+$, the symmetric Wilson polynomial $E^+(\cdot,\ga_{2n})$ is a solution to the difference equation $(Y^2)_{sym} f = \ga_{2n}^2 f$. Now write
\[
\ga_{2n}^2-(t_0+t_1)^2= (t_0+t_1+n)^2-(t_0+t_1)^2=n(2t_0+2t_1+n).
\]
to see that $E^+(x,\ga_{2n})$ is an eigenfunction of $L=(Y^2)_{sym}-(t_0+t_1)^2$ for eigenvalue $n(n+a+b+c+d-1)$. 

Now let us obtain the explicit expression for $(Y^2)_{sym}$. From the explicit expression for $T_0$ and $T_1$,
\[
T_i = t_i + c_{i}(\cdot)(s_i-1),
\] 
and from $s_0 = \tau(-1)s_1= s_1\tau(1)$  we see that we can write 
\[
(Y^2)_{sym} = B(x)[\tau(1)-1] + C(x)[\tau(-1)-1] + D(x),
\]
for $B,C,D \in \C(x)$. To find $D$ it is enough to calculate $Y^2 1$. From $T_i 1= t_i$, since $1 \in \mathcal A^{W}$, we obtain $D(x) =  (t_0+t_1)^2$.
To find $B(x)$ we need to find the coefficient of $\tau(1)$ and of $\tau(1)s_1$ in $T_0T_1+T_1T_0$. The only contribution comes from the terms with $\tau(1)=s_1s_0$, so from $T_1T_0$, and we find 
\[
B(x)=c_{0}(-x)c_{1}(x).
\]
To find $C(x)$ we need to find the coefficients of $\tau(-1)$ and $\tau(-1)s_1$ in $T_0T_1+T_1T_0$. Then we obtain 
\[
C(x)= c_{0}(x)\big(2t_1-c_{1}(x)\big) = c_{0}(x)c_{1}(-x). 
\]

{\rm (b)} Using the duality property for the symmetric Wilson polynomials we find from {\rm (a)} 
\[
\begin{split}
&\big((x_{2m}^2-(t_1+u_1)^2 \big)  E^+(x_{2m}, \ga_{2n})=\\ 
&A^\si(\ga_{2n}) \Big(E^+(x_{2m},\ga_{2n}+1)- E^+(x_{2m},\ga_{2n}) \Big) 
+ A^\si(-\ga_{2n})\Big(E^+(x_{2m},\ga_{2n}-1)- E^+(x_{2m},\ga_{2n}) \Big).
\end{split}
\]
This identity holds for all $m \in \Z_{\geq 0}$, so it holds as an identity in $\mathcal A$. Setting $B_n=A^\si(\ga_{2n})$, $C_n=A^\si(-\ga_{2n})$, and using $\ga_{2n} \pm 1 = \ga_{2n\pm 2}$, gives the desired relation.
\end{proof}
Proposition \ref{prop:diff eq} gives precisely the difference equation and recurrence relation for the well-known Wilson polynomials \cite{Wi80}. This gives us the following expression.
\begin{thm}
The symmetric Wilson polynomials have the explicit expression
\[
E^+(x,\ga_{2n}) = \F{4}{3}{-n, n+a+b+c+d-1, a+x, a-x}{a+b, a+c, a+d}{1}.
\]
\end{thm}
Here we use the standard notation for hypergeometric series,
\[
\F{p}{q}{\al_1, \ldots, \al_p}{\be_1, \ldots, \be_q}{x} = \sum_{j=0}^\infty \frac{ (\al_1)_j \cdots (\al_p)_j }{ (\be_1)_j \cdots (\be_q)_j } \frac{ x^j}{j!}.
\]
Using $E^+(x,\ga_{2n})=P_{2n}^+(x)/P_{2n}^+(x_0)$ and the evaluation of $P^+_{2n}(x_0)$ from Proposition \ref{prop:evaluation}, we now also have an explicit expression for $P_{2n}^+(x)$. Then from Lemma \ref{lem:sym pol} and Theorem \ref{thm:shift} we find an explicit expression for the non-symmetric Wilson polynomial $p_m$ as a sum of two balanced $_4F_3$-series.

\section{The polynomial Fourier transform} \label{sec4}
Let $V$ and $\widetilde V$ be a $\mathcal H$-module and a $\mathcal H_\si$-module, respectively. A Fourier transform associated to the duality isomorphism $\si$ is a linear map $\mathbb F : V \rightarrow \widetilde V$ that intertwines the actions of $\mathcal H$ and $\mathcal H_\si$, i.e.,
\begin{equation} \label{eq:Ftransform}
\mathbb F \circ X = \si(X) \circ \mathbb F.
\end{equation}
We are interested in Fourier transforms that can be written as integral transforms with some kernel. In this section we consider a Fourier transform with the non-symmetric Wilson polynomials as a kernel. In the next section we consider a Fourier transform with a non-polynomial kernel.

\subsection{The Fourier transform $\mathbb F$}
Let $\Ga$ denote the spectrum of $-Y \in \mathcal H$, i.e., $\Ga=\{-\ga_m \ | \ m \in \Z_{\geq 0} \}$, and let $F$ be the space of complex functions on $\Ga$ with finite support,
\[
F = \Big\{ f: \Ga \rightarrow \C \ \vert \ \mathrm{supp}(f) \text{ finite} \Big\}.
\]
We define the non-symmetric polynomial Fourier transform $\mathbb F=\mathbb F_{\mathbf t}:\mathcal A \rightarrow F$ by
\[
(\mathbb Fp)(\ga) = \langle p, E(\cdot,-\ga;\mathbf{t}) \rangle_{\mathbf t}, \qquad \ga \in \Ga.
\]
The mapping $\mathbb F$ is injective, since the bilinear form $\langle\cdot,\cdot \rangle_{\mathbf t}$ is non-degenerate, and, since the polynomials $E(\cdot,\ga_m;\mathbf t)$ form an orthogonal basis for $\mathcal A$ with respect to $\langle\cdot,\cdot \rangle_{\mathbf t}$, $\mathbb F$ is also surjective. For any $f \in F$ let us define the values of $f$ at $\ga_0$ by $f(\ga_0)=f(-\ga_0)$. Now we define an action of the affine Weyl group $\mathcal W$ on $F$ in the same way as on $\mathcal A$, i.e.,
\[
(s_0 f)(\ga) = f(1-\ga), \qquad (s_1f)(\ga)=f(-\ga), \qquad \ga \in \Ga.
\]
\begin{prop} \label{prop:A<->F}
The applications
\[
\begin{split}
(T_i^\si f)(\ga) &= t_i^\si f(\ga) + c_i(\ga;\mathbf t^\si)\big( (s_if)(\ga) - f(\ga) \big), \qquad i=0,1,\\
(p(z)f)(\ga) &= p(\ga)f(\ga), \qquad p \in \mathcal A,
\end{split}
\]
extend uniquely to a representation of $\mathcal H_\si$ on $F$. With this representation, $\mathbb F:\mathcal A \rightarrow F$ is a Fourier transform associated to $\si$.
\end{prop}
\begin{proof}
Assuming that the intertwining property \eqref{eq:Ftransform} holds for $X_1, X_2 \in \mathcal H$, we have for $p \in \mathcal A$,
\[
\mathbb F(X_1 X_2 p) = \si(X_1) \big(\mathbb F (X_2 p) \big) = \si (X_1) \si(X_2) (\mathbb F p) = \si(X_1 X_2)(\mathbb F p),
\]
since $\si$ is an algebra isomorphism. So it is enough to check the intertwining property for $Y, T_1, U_1 \in \mathcal H$, since these elements generate $\mathcal H$ as an algebra.

Let $p \in \mathcal A$. Then for $\ga \in \Ga$,
\[
\mathbb F(Yp)(\ga) = \langle Yp, E(\cdot,-\ga) \rangle_{\mathbf t} =  \langle p, YE(\cdot,-\ga) \rangle_{\mathbf t} = -\ga (\mathbb F p)(\ga).
\]
So we obtain
\[
\mathbb F(Yp) = -z(\mathbb F p) = \si(Y)(\mathbb F p).
\]
For the action of $T_1$ we use Lemma \ref{lem:h1},
\[
\begin{split}
\big(\mathbb F( T_1p) \big)(\ga) &= \langle T_1p, E(\cdot,-\ga) \rangle_{\mathbf t} \\
&= t_1 \langle p, E(\cdot,-\ga) \rangle_{\mathbf t} + c_{1}^\si(\ga) \langle p , E(\cdot,\ga) - E(\cdot,-\ga) \rangle_{\mathbf t}\\
&=  t_1 (\mathbb F p)(\ga) + c_{1}^\si(\ga)\big((\mathbb F p)(-\ga) - (\mathbb F p)(\ga)\big) \\
&= \big(T_1^\si (\mathbb F p)\big)(\ga),
\end{split}
\]
so $\mathbb F(T_1p) = \si(T_1)(\mathbb Fp)$. Recall here that $c^\si(-\ga_0)=0$ and $E(x,-\ga_0)=1$. 

Finally, for $U_1$ we have by Lemma \ref{lem:h1}
\[
\begin{split}
\big(\mathbb F( U_1 p) \big)(\ga) &= \langle U_1 p, E(\cdot,-\ga) \rangle_{\mathbf t} \\
&= u_1 \langle p, E(\cdot,-\ga) \rangle_{\mathbf t} + c_{0}^\si(\ga) \langle p_2 , E(\cdot,-1+\ga) - E(\cdot,-\ga) \rangle_{\mathbf t}\\
&= u_1 (\mathbb F p)(\ga) + c_{0}^\si(\ga)\big((\mathbb F p)(1-\ga) - (\mathbb F p)(\ga)\big) \\
&= \big(T_0^\si (\mathbb F p)\big)(\ga),
\end{split}
\]
which is the same as $\mathbb F(U_1p) = \si(U_1)(\mathbb Fp)$.
\end{proof}

\subsection{The inverse transform}
We define a weight function $w$ on $\Ga$ by
\[
w(\ga;\mathbf{t}) = \Res{y=\ga} \De(y;\mathbf{t}^\si).
\]
Let us associate parameters $\tilde a, \tilde b, \tilde c, \tilde d$ to the multiplicity function $\mathbf t^\si$ in the same way as \eqref{eq:abcd}, i.e.,
\[
(\tilde a, \tilde b, \tilde c, \tilde d)=(t_1+t_0, t_1-t_0, u_1+u_0+\hf, u_1-u_0+\hf).
\]
In view of the explicit expression for $\De$, $\Res{y=-n}\Ga(y)=(-1)^n/n!$, $n \in \Z_{\geq 0}$, and $\Ga(n+1)=n!$, we have
\[
\begin{split}
w&(\ga;\mathbf{t})=\\
& \begin{cases} 
\displaystyle 
(-1)^{n-1}\frac{ \Ga(\tilde a+\ga) \Ga(\tilde b+\ga) \Ga(\tilde b +1-\ga) \Ga(\tilde c \pm \ga) \Ga(\tilde d \pm \ga)} { \Ga(\ga - \tilde a) \Ga(2\ga) \Ga(1-2\ga) }, & \quad \ga = \tilde a+n,\quad n \in \N,\\ \\
\displaystyle 
(-1)^n \frac{ \Ga(\tilde a +1 - \ga) \Ga(\tilde b+\ga) \Ga(\tilde b +1 -\ga) \Ga(\tilde c \pm \ga) \Ga( \tilde d\pm \ga) }{ \Ga(1- \ga-\tilde a ) \Ga(2 \ga) \Ga(1-2\ga) }, & \quad \ga=-(\tilde a +n ), \quad n \in \Z_{\geq 0}.
\end{cases}
\end{split}
\]
To this weight function we associate a bilinear form $[\cdot,\cdot]_{\mathbf{t}}: F \times F \rightarrow \C$ by
\[
[f,g]_{\mathbf{t}} = \sum_{\ga \in \Ga} f(\ga)g(-\ga) w(\ga;\mathbf t).
\]
Now we define the map $\mathbb G=\mathbb G_{\mathbf{t}} : F \rightarrow \mathcal A$ by
\[
(\mathbb G f)(x)  = [f,E(x,\cdot;\mathbf t)]_{\mathbf t}, \qquad f \in F, \ x \in \C.
\]
Note here that $E(x,\cdot) \not\in F$, but clearly for any $f \in F$ the function $\mathbb G f$ exists.
\begin{prop} \label{prop:F<->A}
For $X \in \mathcal H_\si$ and $f \in F$ we have
\[
\mathbb G(X f) = \si^{-1}(X) (\mathbb Gf).
\]
\end{prop}
\begin{proof}
Let $f \in F$. It is enough to check the statement for $X=z, T_0^\si, T_1^\si$. For $X=z$ we obtain
\[
\begin{split}
\mathbb G(z f)(x) =& \sum_{\ga \in \Ga} \ga f(\ga) E(x,-\ga) w(\ga)  \\
=& \sum_{\ga \in \Ga} f(\ga) \big( -Y E(\cdot,-\ga)\big)(x) w(\ga) \\
=& \big( -Y(\mathbb G f)\big)(x)\\
=& \big(\si^{-1}(z) (\mathbb G f) \big)(\ga).
\end{split}
\]
Next, for $X=T_0^\si$ we obtain
\[
\begin{split}
\mathbb G&(T_0^\si f)(x)\\
= & \sum_{\ga \in \Ga} \Big(u_1 f(\ga) + \frac{ (u_1 + u_0 + \hf-\ga)(u_1-u_0+\hf-\ga)}{1-2\ga}\big(f(1-\ga)-f(\ga)\big) \Big) E(x,-\ga) w(\ga) \\
=& \sum_{\ga \in \Ga} f(\ga) \Big(u_1 E(x,-\ga) + \frac{ (\tilde c-\ga)(\tilde d-\ga)}{1-2\ga} \big(E(x,\ga-1) - E(x,-\ga)\big) \Big) w(\ga) \\
=& \sum_{\ga \in \Ga} f(\ga) \big(U_0 E(\cdot;-\ga) \big)(x) w(\ga) \\
=& \big(\si^{-1}(T_0^\si) (\mathbb G f) \big)(x)
\end{split}
\]
Here we use Lemma \ref{lem:h1}, and 
\[
\frac{ (\tilde c-\ga)(\tilde d-\ga)}{1-2\ga}w(\ga;\mathbf{t}) = \frac{ (\tilde c-1+\ga)(\tilde d-1+\ga)}{2\ga-1}w(1-\ga;\mathbf{t}),
\]
which can be obtained from the explicit expression for $w(\ga;\mathbf{t})$. 

Finally, using Lemma \ref{lem:h1} we find for $T_1^\si$,
\[
\begin{split}
\mathbb G&(T_1^\si f)(x) \\
=& \sum_{ \ga \in \Ga} \Big(t_1 f(\ga) + \frac{ (t_1+t_0 + \ga)(t_1-t_0+\ga)}{2\ga} \big( f(-\ga) - f(\ga)  \big)  \Big) E(x,-\ga) w(\ga)\\
=& \sum_{ \ga \in \Ga}f(\ga) \Big( t_1 E(x,-\ga) + \frac{ (\tilde a + \ga) (\tilde b+ \ga) }{2\ga}\big( E(x, \ga) - E(x, -\ga) \big) \Big) w(\ga)\\
=& \sum_{ \ga \in \Ga}f(\ga)  \big(T_1 E(\cdot,-\ga)(x) w(\ga)\\
=& \big(\si^{-1}(T_1^\si) (\mathbb G f) \big)(x),
\end{split}
\]
where we have used
\[
\frac{ (\tilde a + \ga) (\tilde b+ \ga) }{2\ga} w(\ga;\mathbf t) = \frac{ (\tilde a - \ga) (\tilde b- \ga) }{-2\ga} w(-\ga;\mathbf t),
\]
for $\ga\neq -\ga_0$.
This proves the proposition.
\end{proof}

From combining Propositions \ref{prop:A<->F} and \ref{prop:F<->A} it follows that $\mathbb G$ is, up to a constant, the inverse of $\mathbb F$. This enables us to determine the `quadratic norms' for the non-symmetric Wilson polynomials.
\begin{thm} \label{thm:FG}
{\rm (a)} We have $\mathbb G \circ \mathbb F = \mathcal N \, id_{\mathcal A}$, and $\mathbb F \circ \mathbb G = \mathcal N \, id_{F}$, where $\mathcal N =\mathcal N(\mathbf t)= \langle 1,1\rangle_{\mathbf t}\, w(-\ga_0;\mathbf{t})$.

{\rm (b)} For $\ga \in -\Ga$ we have
\[
\frac{\langle E(\cdot,\ga;\mathbf t), E(\cdot,\ga; \mathbf t) \rangle_{\mathbf t}}{\langle 1, 1 \rangle_{\mathbf t}} = \frac{ w(-\ga_0;\mathbf t) }{ w(-\ga;\mathbf t) }.
\]
\end{thm}
\begin{proof}
{\rm (a)} Let $p \in \mathcal A$. From Propositions \ref{prop:A<->F} and \ref{prop:F<->A} we obtain
\[
\mathbb G(\mathbb F p) = \mathbb G\big( \si(p(z))(\mathbb F 1) \big) = p(z) (\mathbb G (\mathbb F 1)).
\]
Using the orthogonality relations for the non-symmetric Wilson polynomials, and the definitions of $\mathbb F$ and $\mathbb G$, we obtain
\[
\mathbb G(\mathbb F 1) = \langle 1, 1 \rangle_{\mathbf t} w(-\ga_0;\mathbf t).
\]
So we have $\mathbb G \circ \mathbb F = \mathcal N \, id_{\mathcal A}$ with $\mathcal N$ as in the theorem. 

Next let $f \in F$. Since $\mathbb F$ is bijective, we have $f =\mathbb Fp$ for some $p \in \mathcal A$. This gives us
\[
\mathbb F (\mathbb G f) = (\mathbb F \mathbb G)(\mathbb F p) = \mathcal N (\mathbb F p) = \mathcal N f,
\]
which shows that $\mathbb F \circ \mathbb G = \mathcal N \, id_{F}$.

{\rm (b)} Using the orthogonality relations for the non-symmetric Wilson polynomials,
and the definitions of $\mathbb F$ and $\mathbb G$, we find
\[
\big(\mathbb G(\mathbb F E(\cdot,-\ga))\big)(x) = \langle E(\cdot,-\ga),E(\cdot,-\ga) \rangle_{\mathbf t} E(x;-\ga)w(\ga;\mathbf t), \qquad \ga \in \Ga.
\]
On the other hand, from {\rm (a)} we obtain $\big(\mathbb G(\mathbb F E(\cdot,-\ga))\big)(x) = \mathcal N E(x,-\ga)$. Comparing gives
\[
\langle E(\cdot,-\ga),E(\cdot,-\ga) \rangle_{\mathbf t} = \frac{ \mathcal N }{ w(\ga;\mathbf t)},
\]
which is equivalent to statement {\rm (b)} in the theorem.
\end{proof}

Let us remark that since $\langle f,g \rangle_{\mathbf t} = \hf(a+b) \langle f,g \rangle_{\mathbf t}^+$ for $f,g \in \mathcal A^W$, and $\langle 1, 1 \rangle_{\mathbf t}^+$ is a well-known integral (see \cite{Wi80}), we can evaluate $\langle 1, 1 \rangle_{\mathbf t}$;
\begin{equation} \label{eq:<1,1>}
\langle 1, 1 \rangle_{\mathbf t} = \frac{ \Ga(a+b+1) \Ga(a+c) \Ga(a+d) \Ga(b+c) \Ga(b+d) \Ga(c+d) }{ \Ga(a+b+c+d) }.
\end{equation}
So the `quadratic norms' of the non-symmetric Wilson polynomials are completely explicit.

For $\mathbb F$ and $\mathbb G$ we have the following Plancherel-type formulas.
\begin{prop}
For $p_1,p_2 \in \mathcal A$ and $f_1,f_2 \in F$,
\[
[\mathbb F p_1, \mathbb F p_2]_{\mathbf t} = \mathcal N \langle p_1, p_2 \rangle_{\mathbf t}, \qquad \langle \mathbb Gf_1, \mathbb Gf_2 \rangle_{\mathbf t} = \mathcal N [f_1, f_2 ]_{\mathbf t}.
\]

\end{prop}
\begin{proof}
From the definitions of $\mathbb F$ and $\mathbb G$ we obtain, for $p \in \mathcal A$ and $f \in F$,
\begin{equation} \label{eq:planch}
\begin{split}
[f, \mathbb Fp]_{\mathbf t} &= \sum_{\ga \in \Ga} f(\ga) \left(  \frac{1}{2\pi i} \int_{\mathcal C} p(x) E(x, -\ga) \De(x) dx \right) w(\ga)\\
&= \frac{1}{2\pi i} \int_{\mathcal C} \left( \sum_{\ga \in \Ga} f(\ga) E(x,-\ga) w(\ga) \right) p(x) \De(x) dx \\
&=\langle \mathbb Gf , p \rangle_{\mathbf t}. 
\end{split}
\end{equation}
For the first identity take $f= \mathbb F p_1$ and $p=p_2$ in \eqref{eq:planch}, for the second identity take $f=f_1$ and $p = \mathbb G f_2$. Then use Theorem \ref{thm:FG}. 
\end{proof}

\subsection{The symmetric transform}
The non-symmetric polynomial Fourier transform is defined as an integral transform with the non-symmetric Wilson polynomial as a kernel. Similarly we can define a symmetric polynomial Fourier transform as an integral transform with the symmetric Wilson polynomial as a kernel. For this transform we can also find the inverse, and obtain Plancherel-type formulas. The proofs for these are very similar to the proofs for the non-symmetric case, and will be left to the reader.

Let $\Ga^+ = \{ \ga_{2n}\ | \ n \in \Z_{\geq 0}\} \subset -\Ga$. Furthermore, let $F_+$ denote the subspace of $F$ consisting of $W$-invariant functions on $\Ga$, then
\[
F_+ = \{ f \in F \ | \ f(\ga)=f(-\ga), \ \ga \in \Ga^+\}.
\]
We define the symmetric polynomial Fourier transform $\mathbb F^+: \mathcal A_+ \rightarrow  F_+$ by
\[
(\mathbb F^+p)(\ga) = \langle p, E^+(\cdot,\ga)\rangle_{\mathbf t}^+, \qquad \ga \in \Ga^+, \ p \in \mathcal A_+.
\] 
Let $w^+$ be the weight function on $\Ga^+$ defined by
\[
w^+(\ga;\mathbf t) = \Res{y=\ga} \De^+(y;\mathbf t^\si).
\]
With this weight function we define a bilinear form $[\cdot,\cdot]_{\mathbf t}^+: F_+ \times F_+ \rightarrow \C$ by
\[
[f,g]_{\mathbf t}^+ = \sum_{ \ga \in \Ga^+} f(\ga) g(\ga) w^+(\ga;\mathbf t).
\]
Now we define $\mathbb G^+ : F_+ \rightarrow \mathcal A_+$ by
\[
(\mathbb G^+ f)(x) = [f, E^+(x,\cdot)]^+_{\mathbf t}, \qquad x \in \C,\ f \in F_+.
\]
In the same way as in the previous section we can show that $\mathbb G^+$ is, up to a multiplicative constant, the inverse of $\mathbb F^+$, and this leads to the (well-known) evaluation of the quadratic norms of the symmetric Wilson polynomials and Plancherel-type formulas. 
\begin{thm} \label{thm:symm}
{\rm (a)} We have $\mathbb G^+ \circ \mathbb F^+ = \mathcal N^+ id_{\mathcal A_+}$, and $\mathbb F^+ \circ \mathbb G^+ = \mathcal N^+ id_{F_+}$, where $\mathcal N^+ = \mathcal N^+(\mathbf t) = \langle 1, 1 \rangle_{\mathbf t}^+ w^+(\ga_0;\mathbf t)$.

{\rm (b)} For $\ga \in \Ga^+$,
\[
\frac{ \langle E^+(\cdot, \ga ;\mathbf t),E^+(\cdot, \ga ;\mathbf t) \rangle_{\mathbf t}^+ }{ \langle 1, 1 \rangle_{\mathbf t}^+ } = \frac{ w^+(\ga_0;\mathbf t)}{ w^+(\ga;\mathbf t)}.
\]

{\rm (c)} For $p_1, p_2 \in \mathcal A_+$ and $f_1, f_2 \in F_+$,
\[
[ \mathbb F^+p_1, \mathbb F^+p_2]_{\mathbf t}^+ = \mathcal N^+ \langle p_1, p_2 \rangle_{\mathbf t}^+, \qquad \langle \mathbb G^+ f_1, \mathbb G^+ f_2 \rangle_{\mathbf t}^+ = \mathcal N^+ [ f_1, f_2]_{\mathbf t}^+.
\] 
\end{thm}

\section{The non-polynomial Fourier transform} \label{sec5}
In this section we study a Fourier transform associated to $\si$ with a non-polynomial kernel. This kernel turns out to be a meromorphic continuation of the non-symmetric Wilson polynomial in its degree.

\subsection{The Gaussian}
We define an involution $\tau$ on multiplicity functions by interchanging the values on the $a_0$-orbit and the $a_0^\vee$-orbit; given $\mathbf t=(t_0,u_0,t_1,u_1)$ we have
\[
\mathbf t^{\tau} = (u_0, t_0, t_1, u_1).
\]
We use notations with $\tau$ in the same way we do with the involution $\si$. We will also compose $\si$ and $\tau$, and we use notations like $\mathbf t^{\si \tau} = (\mathbf t^\si)^\tau=(u_0, u_1, t_1, t_0)$. We will use similar notations for objects depending on $\mathbf t^{\si \tau}$, e.g.~$H_{\si \tau} =  H(\mathbf t^{\si \tau})$. We mention that
\[
\mathbf t^{ \tau \si \tau} = (t_0,u_1,t_1,u_0)= \mathbf t^{\si \tau \si},
\]
which will be used later on.

We define the Gaussian by
\[
G(x; \mathbf t) = \Ga(t_0 - u_0 + \hf  \pm x) .
\]
This function will play the same role as $e^{-x^2}$ does for the Hankel transform, see \cite{CM02}. Let $\tau$ denote conjugation with the Gaussian;
\begin{equation} \label{def:tau}
G \circ X = \tau(X) \circ G, \qquad X \in \mathcal H.
\end{equation}
\begin{prop}
The application $X \mapsto \tau(X)$ is an algebra isomorphism $\mathcal H \rightarrow \mathcal H_\tau$. On the generators $z,T_0, T_1$ of $\mathcal H$ the action of $\tau$ is given by
\[
\tau(z) = z, \quad \tau(T_0) = U_0^\tau, \quad \tau(T_1)=T_1^\tau.
\]
\end{prop}
\begin{proof}
The fact that $\tau$ is an algebra isomorphism follows from the explicit action on the generators, so we need to check the actions given in the lemma.

The identities $GzG^{-1}=z$ and $GT_1G^{-1}=T_1^\tau$ are easy to check, using that $z$ and $T_1$ do not depend on $t_0$ and $u_0$, and the identity $(s_1G^{-1})(x)=G(x)^{-1}$. For $T_0$ we use
\[
(s_0 G^{-1})(x) = G(1-x)^{-1} = \frac{ t_0-u_0-\hf +x }{ t_0 - u_0 +\hf -x } G(x)^{-1},
\]
then we obtain, for $p \in \mathcal A$,
\[
\big(\tau(T_0) p\big)(x) = \big(t_0 - c_{0}(x)\big)p(x) + \frac{ (t_0+u_0+\hf -x)(t_0-u_0-\hf +x) }{1-2x}  p(1-x).
\]
By a direct calculation it is verified that
\[
\begin{split}
t_0-c_{0}(x) = x-\hf-u_0 + c_0^\tau(x).
\end{split}
\]
So we find
\[
\big(\tau(T_0) p\big)(x) = (x-\hf-u_0)p(x) - c_{0}^\tau(x)\big(p(1-x)-p(x) \big) = \big((z-\hf-T_0^\tau)p\big)(x). \qedhere
\]
\end{proof}

\subsection{Construction of a Fourier transform}
We construct a Fourier transform associated to $\si$ using the transforms $\mathbb F$, $\mathbb G$ from the previous section, and the Gaussian. First, we define a linear operator $\mathfrak F: \mathcal A \rightarrow \mathcal A$ by
\[
\mathfrak F = \mathbb G_{\si \tau} \circ G_{\tau \si} \circ \mathbb F_{\tau}.
\]
Note that the set $\Ga_\tau$ consists of the points
\[
-\ga_m^\tau = \begin{cases}
-(u_0+t_1 +n), & m=2n,\\
u_0+t_1 + n, & m=2n-1,
\end{cases}
\]
which are clearly invariant under $u_1 \leftrightarrow t_0$. So we have $\Ga_\tau = \Ga_{\si \tau}$. Recall that $F$ is the set of finitely supported functions on $\Ga$, then it follows that $F^\tau = F^{\si \tau}$. So the composition $\mathbb G_{\si \tau} \circ G_{\tau \si} \circ \mathbb F_{\tau}$ is well-defined. 
\begin{prop} \label{prop:properties1}
For $X \in \mathcal H_\tau$ we have
\[
\mathfrak F \circ X = \chi(X) \circ \mathfrak F,
\]
where $\chi:\mathcal H_\tau \rightarrow \mathcal H_{\si \tau}$ is the isomorphism $\chi = \si^{-1}_{\si \tau} \circ \tau_{\tau \si} \circ \si_{\tau}$. Explicitly, $\chi : \mathcal H_\tau \rightarrow \mathcal H_{\si \tau}$ is given on algebraic generators of $\mathcal H_\tau$ by
\[
\chi(z) = -U_0^{\si \tau}- T_1^{\si\tau}, \qquad \chi(Y^\tau)=Y^{\si \tau}, \qquad \chi(T_1^\tau) = T_1^{\si \tau}.
\]
\end{prop}
\begin{proof}
First we check the intertwining property. Let $X \in \mathcal H_\tau$ and $p \in \mathcal A$. Using the intertwining properties of $\mathbb F$ and $G$, see Proposition \ref{prop:A<->F} and \eqref{def:tau}, we have
\[
\begin{split}
\mathfrak F (Xp) &= (\mathbb G_{\si \tau} \circ G_{\tau \si} \circ \mathbb F_\tau)(Xp) \\
&= \big(\mathbb G_{\si \tau} \circ G_{\tau \si}\big)\big( \si_{\tau}(X) (\mathbb F_\tau p) \big) \\
&= \mathbb G_{\si \tau} \Big[ (\tau_{\tau \si} \circ \si_{\tau})(X) \big( (G_{\tau \si} \circ \mathbb F_\tau)(p) \big)\Big].
\end{split}
\]
Note that $(\tau_{\tau \si} \circ \si_{\tau})(X) \in \mathcal H_{\tau \si \tau} = \mathcal H_{\si \tau \si}$. So by the intertwining property of $\mathbb G$, see Proposition \ref{prop:F<->A}, it follows that
\[
\mathfrak F(Xp) =  (\si_{\si \tau}^{-1} \circ \tau_{\tau \si} \circ \si_{\tau})(X) \big((\mathbb G_{\si \tau} \circ G_{\tau \si} \circ \mathbb F_\tau)(p) \big) = \chi(X) (\mathfrak F p). 
\]
The explicit expressions for $\chi(X)$, $X=z,Y^\tau, T_1^\tau \in \mathcal H_\tau$, are easily verified using the explicit actions of  $\si$ and $\tau$ on algebraic generators. 
\end{proof}
We define the symmetric version $\mathfrak F^+: \mathcal A_+ \rightarrow \mathcal A_+$ of the operator $\mathfrak F$ by
\[
\mathfrak F^+ = \mathbb G_{\si \tau}^+ \circ G_{\tau \si} \circ \mathbb F_{\tau}^+.
\]
\begin{prop} \label{prop:properties2} 
{\rm (a)} Let $\mathfrak L:\mathcal A \rightarrow \mathcal A$ be the linear operator defined on the Wilson polynomials by
\[
\big(\mathfrak L E_\tau(\cdot,\ga)\big)(\la) =  G_{\tau \si \tau}(\ga)  E_{\si \tau}(\la,\ga), \qquad \ga \in -\Ga_\tau, 
\]
then $\mathfrak L= \mathcal N_\tau^{-1} \, \mathfrak F$, where $\mathcal N$ denotes the constant from Theorem \ref{thm:FG}.

{\rm (b)} Let $\mathfrak L^+:\mathcal A_+ \rightarrow \mathcal A_+$ be the linear operator defined on the symmetric Wilson polynomials by
\[
\big(\mathfrak L^+ E_\tau^+(\cdot,\ga)\big)(\la) = G_{\tau \si \tau}(\ga)  E_{\si \tau}^+(\la,\ga), \qquad \ga \in \Ga_\tau^+,
\]
then $\mathfrak L^+= (\mathcal N_{\tau}^+)^{-1} \, \mathfrak F^+$, where $\mathcal N^+$ denotes the constant from Theorem \ref{thm:symm}.
\end{prop}
\begin{proof}
{\rm (a)} Let $\mathcal L$ and $\mathcal L'$ be linear operators on $\mathcal A$, such that $\mathfrak LE_\tau(\cdot,\ga)=\mathfrak L' E_\tau(\cdot,\ga)$. The Wilson polynomials $E_\tau(\cdot,\ga)$, $\ga\in -\Ga_\tau$, form a linear basis for $\mathcal A$, therefore $\mathfrak L=\mathfrak L'$. So, to prove the proposition, we need to calculate the action of $\mathfrak F$ on the Wilson polynomials $E_\tau$.

Let $\ga, \ga' \in \Ga_\tau$. By Theorem \ref{thm:FG}{\rm (b)} we have
\[
\big(\mathbb F_\tau E_\tau(\cdot,-\ga)\big)(\ga') = \langle E_\tau(\cdot,-\ga), E_\tau(\cdot,-\ga') \rangle_{\mathbf t^\tau} = \de_{\ga,\ga'}  \frac{ \mathcal N_{\tau} }{w_\tau(\ga)}.
\]
Multiplying this expression by $G_{\tau \si}(\ga')$ and applying $\mathbb G_{\si \tau}$ then leads to
\[
\big(\mathfrak F E_\tau(\cdot,-\ga)\big)(\la) = \mathcal N_{\tau}  G_{\tau \si}(\ga) \frac{w_{\si \tau}(\ga)}{w_{\tau}(\ga)} E_{\si \tau}(\la,-\ga).
\]
Using the definition of the weight $\De$ we see that $G_\tau \De = G\De_\tau$, 
so that by the definition of the weight $w$ and by $\mathbf t^{\tau \si}= \mathbf t^{\si \tau \si \tau}$ it follows that
\begin{equation} \label{eq:Gw}
G_{\tau \si}(\ga) w_{\si \tau}(\ga) = G_{\tau \si \tau}(\ga) w_\tau(\ga), \qquad \ga \in \Ga_\tau.
\end{equation}
So we find that $\mathfrak F E_\tau(\cdot,-\ga)=\mathcal N_\tau  \mathfrak L E_\tau(\cdot, -\ga)$.

{\rm (b)} This is proved in the same way as {\rm (a)} using the orthogonality relations for the symmetric Wilson polynomials.
\end{proof}

Now define the operator $\mathcal F' : \mathcal AG_\tau \rightarrow \mathcal AG_{\si\tau}$ by
\begin{equation} \label{eq:def F}
\mathcal F'  = \,G_{\si \tau}\circ \mathfrak F \circ G_{\tau}^{-1}.
\end{equation}
\begin{prop} \label{prop:F}
The operator $\mathcal F'$ is a Fourier transform associated to $\si$.
\end{prop}
\begin{proof}
We need to show that $\mathcal F' \circ X = \si(X) \circ \mathcal F'$ for all $X \in \mathcal H$. Let $p \in \mathcal A$ and let $f=pG_{\tau}$. From \eqref{def:tau} and Proposition \ref{prop:properties1}{\rm (a)} we obtain
\[
\begin{split}
\mathcal F'(Xf) &= G_{\si \tau} \mathfrak F'(\tau(X)p)=  G_{\si\tau} [(\chi\circ \tau)(X)] (\mathfrak Fp)\\
& = [(\tau_{\si\tau}\circ \chi\circ \tau)(X)] \big(G_{\si\tau} \mathfrak Fp\big) =k\, [(\tau_{\si\tau}\circ \chi\circ \tau)(X)] (\mathcal F'f).
\end{split}
\]
Now we must check that $\si= \tau_{\si\tau}\circ \chi\circ \tau$. It is enough to check that this identity is true on generators of $\mathcal H$. This is a straightforward exercise that we leave to the reader.
\end{proof}

\subsection{The Fourier transform $\mathcal F$} 
We want to write $\mathcal F'$ as an integral transform with some kernel. First we consider the related transform $\mathfrak F$. Writing out explicitly $\mathfrak Fp= (\mathbb G_{\si \tau} \circ G_{\tau \si} \circ \mathbb F_\tau)p$ and changing the order of summation and integration, suggests that we may write
\begin{equation} \label{eq:inttrans}
(\mathfrak F p)(\la) = k\langle p, \mathfrak E(\cdot,\la) \rangle_{\mathbf t^{\tau}}, \qquad p \in \mathcal A,
\end{equation}
for some kernel $\mathfrak E(x,\la)$ and $k$ a non-zero constant.
\begin{rem}
If we formally expand the kernel $\mathfrak E(x,\la)$ in terms of Wilson polynomials $E_\tau(x,\ga)$ and we calculate the coefficients of $E_\tau(x,\ga)$ using the orthogonality relations for the Wilson polynomials, we obtain
\[
\mathfrak E(x,\la) = k^{-1} \sum_{\ga \in \Ga_\tau} E_\tau(x; -\ga) E_{\si \tau}(\la;-\ga) G_{\tau \si}(\ga) w_{\si \tau}(\ga).
\]
This gives a formal power series for the kernel $\mathfrak E(x,\la)$ similar to the kernels defined by Cherednik \cite[(5.12)]{Ch97}, see also \cite[(5.13)]{St03}. However, it turns out here that this sum does not converge absolutely. In \cite{Ch97}, \cite{St03} the absolute convergence of the series comes from the Gaussian, which, in the rank 1 case, contains the factor $q^{n^2}$ for $0<q<1$, where $n$ is the summation index. Although we cannot define the kernel in this way, it does give us an idea what properties the kernel $\mathfrak E$ is expected to have; if the above expansion would converge absolutely, $\mathfrak E(x,\la)$ would be an entire function in $x$ and $\la$, and $\mathfrak E$ would satisfy the duality property $\mathfrak E(x,\la) = \mathfrak E_\si(\la,x)$, provided that $k_\si=k$.
\end{rem}
We want to find a kernel $\mathfrak E$ such that the integral transform \eqref{eq:inttrans} maps a Wilson polynomial $E_\tau(\cdot,\ga)$ to a multiple of $G_{\tau\si \tau}(\ga)E_{\si \tau}(\cdot,\ga)$ as in Proposition \ref{prop:properties2}{\rm (a)}. To do this we introduce the function 
\begin{equation}\label{def:phi}
\begin{split}
\phi_\la(x;\mathbf t)=&\frac{\Ga(1-a-d)}{\Ga(a+b)\Ga(a+c)\Ga(1-d \pm x) \Ga(1-\tilde d \pm \la) }\\
& \times \F{4}{3}{a+x,a-x,\tilde a+x,\tilde a-x}{a+b,a+c,a+d}{1}\\
+&  \frac{ \Ga(a+d-1) }{ \Ga(1+b-d) \Ga(1+c-d) \Ga(a \pm x) \Ga(\tilde a \pm \la) }\\
& \times  \F{4}{3}{1-d+x,1-d-x,1- \tilde d+ \la, 1-\tilde d- \la}{1+b-d, 1+c-d, 2-a-d}{1}.
\end{split}
\end{equation} 
In \cite{Gr03} the function $\phi_\la$ is called a Wilson function. Here we will use the name `Wilson function' for a different (but closely related) function, see Definition \ref{def:nsWilson function}. Using transformation formulas for hypergeometric functions, the function $\phi_\la$ can be expressed as a multiple of a very-well-poised $_7F_6$-function. In \cite{Gr03} a second-order difference operator $L$ is studied, which in the notation of this paper can be written as $L =\tilde a^2-G_\tau^{-1}\circ(Y^2)_{sym} \circ G_\tau$. From results of Ismail, Letessier, Valent and Wimp \cite{ILVW90} and Masson \cite{Ma91}, who studied the associated Wilson polynomials using contiguous relations for $_7F_6$-series, it follows that the function $\phi_\la$ is an eigenfunction of $L$ for eigenvalue $\tilde a^2-\la^2$. 

The following properties of $\phi_\la$ will be useful for us.
\begin{lem} \label{lem:function phi}
The function $\phi_\la(x)$ has the following properties:

{\rm (a)} $\phi_\la(x)$ is an entire function in $x$ and $\la$.

{\rm (b)} For $p \in \mathcal A$ and $y \in \R$ the integral $\langle p, \phi_\la(\,\cdot+y) \rangle_{\mathbf t^\tau}^+$ converges absolutely.

{\rm (c)} For $\ga \in \Ga^+_\tau$
\[
\langle E^+_\tau(\cdot,\ga), \phi_\la \rangle_{\mathbf t^\tau}^+= \frac{2G_{\tau \si \tau}(\ga)}{G_{\tau \si \tau}(\ga_0^\tau)} E^+_{\si \tau}(\la,\ga).
\]
\end{lem}
Note that for (b) the only conditions on the values of $t_0,u_0,t_1,u_1$ are the conditions from subsection \ref{ssec:orth rel}.
\begin{proof}
Property (a) follows directly from the definition of $\phi_\la(x)$. Using a transformation formula for $_7F_6$-series, it is proved in \cite{Gr03} that for $x \rightarrow \pm \infty$ and $y \in \R$,
\[
\phi_\la(y+ix) \sim |x|^{d-a-b-c}e^{\pi (|x|-iy)}\Big( c(\la) |x|^{-2\la} +c(-\la) |x|^{2\la} \Big),
\]
where $c(\la)$ is independent of $x$ and $y$. The weight function $\De_\tau^+$ has the asymptotic behavior
\begin{equation} \label{eq:asymp De}
\De_\tau^+(x) \sim |x|^{2a+2b+2c-2d-1} e^{-2\pi|x|}, \qquad x \rightarrow \pm i\infty,
\end{equation}
which can be obtained from applying Euler's reflection formula for the $\Ga$-function and using the asymptotic formula
\begin{equation} \label{eq:asGa}
\frac{ \Ga(\al+z) }{ \Ga(\be+z) } \sim z^{\al-\be}, \qquad |z| \rightarrow \infty, \ |\arg(z)|<\pi.
\end{equation}
Now we see that the integral $\langle p, \phi_\la(\,\cdot+y) \rangle_{\mathbf t^\tau}^+$ is absolutely convergent, independent of the values of $a,b,c,d$. Property (c) is Theorem 6.9 in \cite{Gr03}, where it is proved in case the parameters $a,b,c,d$ are such that $\De^+_\tau(x)$ is a positive weight on $i\R$, but the prove remains valid without these conditions.
\end{proof}

From Lemma \ref{lem:function phi}{\rm (c)} and Proposition \ref{prop:properties2}{\rm (b)} it follows immediately that $\mathfrak F^+$ can be written as an integral transform with $\phi_\la$ as a kernel. Next we use $\phi_\la$ to write $\mathfrak F$ as an integral transform. Recall the generalized Weyl character formula from Theorem \ref{thm:shift}. For the renormalized Wilson polynomials this formula looks as follows.
\begin{lem} \label{lem:E-+}
For $m \in \N$ the renormalized anti-symmetric Wilson polynomial satisfies
\[
\begin{split}
E^-(x,\ga_{m};\mathbf t) &= (-1)^m \al\,\de_\si(-\ga_{m}) \de(x) E^+(x,\ga_{m-2};t_0,u_0,t_1+1,u_1),
\end{split}
\]
where 
\[
\al=\al(\mathbf t) = \frac{1}{(a+b)(a+b+1)(a+c)(a+d)}.
\] 
\end{lem}
Here we use the convention $E^+(x,\ga_{-1})=1$.
\begin{proof}
From Lemma \ref{lem:sym pol} and $b_{2n}-t_1 = c^\si_1(-\ga_{2n})$ we find
\[
C_-p_{2n} = - \frac{ c_1^\si(-\ga_{2n})}{2t_1} P^-_{2n}, \qquad C_-p_{2n-1} = \frac{1}{2t_1}P_{2n}^-.
\]
So by Theorem \ref{thm:shift} we have $E^-(x,\ga_m)= \al_m \de(x) E^+(x,\ga_{m-2};t_0,u_0,t_1+1,u_1)$ for $m \in \N$. Recall here that $E^+(x,\ga_{2n})=E^+(x,\ga_{2n-1})$. Normalizing the polynomials then gives
\[
\begin{split}
\al_{2n} &= - \frac{ c_1^\si(-\ga_{2n})\, P_{2n-2}^+(x_0+1;t_0,u_0,t_1+1,u_1)} {2t_1\, p_{2n}(-x_0)},\\
\al_{2n-1} &= \frac{P_{2n-2}^+(x_0+1;t_0,u_0,t_1+1,u_1)}{2t_1 \,  p_{2n-1}(-x_0)}.
\end{split}
\]
Writing this out explicitly using Proposition \ref{prop:evaluation} completes the proof.
\end{proof}

\begin{prop} \label{prop:kernel}
Let $\mathfrak E$ be given by
\[
\mathfrak E(x,\la) = \phi_\la(x;\mathbf t) + \de(x) \de_\si(\la) \phi_\la(x;t_0,u_0,t_1+1,u_1),
\]
then 
\begin{equation} \label{eq:id1}
\langle E_\tau(\cdot,\ga), \mathfrak E(\cdot,\la) \rangle_{\mathbf t^\tau} =2t_1\frac{ G_{\tau \si \tau}(\ga)}{ G_{\tau \si \tau}(\ga_0^\tau)}E_{\si \tau}(\la,\ga), \qquad \ga \in -\Ga_\tau.
\end{equation}
Consequently, $(\mathfrak Fp)(\la) = \mathcal K\,\langle p, \mathfrak E(\cdot,\la) \rangle_{\mathbf t^\tau}$ for $p \in \mathcal A$, where the constant  is explicitly given by $\mathcal K =  \mathcal N_\tau G_{\tau \si \tau}(\ga_0^\tau)/2t_1$.
\end{prop}
\begin{proof}
If \eqref{eq:id1} is valid, then by Proposition \ref{prop:properties2}{\rm (a)}  we have $(\mathfrak Fp)(\la) = \mathcal K \langle p, \mathfrak E(\cdot,\la) \rangle_{\mathbf t^\tau}$ for all $p\in \mathcal A$, with $\mathcal K=\mathcal N_\tau G_{\tau \si \tau}(\ga_0^\tau)/2t_1$. So we need to show that \eqref{eq:id1} holds.

Let $\ga \in -\Ga_\tau$. In order to prove \eqref{eq:id1} we split $E_\tau$ and $\mathfrak E$ in symmetric and anti-symmetric parts, then
\[
\langle E_\tau(\cdot,\ga), \mathfrak E(\cdot,\la) \rangle_{\mathbf t^\tau} = \langle E_\tau^+(\cdot,\ga), \mathfrak E^+(\cdot,\la) \rangle_{\mathbf t^\tau} + \langle E_\tau^-(\cdot,\ga), \mathfrak E^-(\cdot,\la) \rangle_{\mathbf t^\tau},
\]
where we denoted $C_\pm^\tau \mathfrak E(\cdot,\la) = \mathfrak E^\pm(\cdot,\la)$. Since $\mathfrak E^+(x,\la)= \phi_\la(x)$, we obtain from Lemmas \ref{lem:blforms} and \ref{lem:function phi}{\rm (c)} 
\[
\langle E_\tau^+(\cdot,\ga), \mathfrak E^+(\cdot,\la) \rangle_{\mathbf t^\tau}
= t_1 \langle E_\tau^+(\cdot,\ga), \mathfrak E^+(\cdot,\la) \rangle_{\mathbf t^\tau}^+ = 2t_1\frac{ G_{\tau \si \tau}(\ga)}{G_{\tau \si \tau}(\ga_0)} E^+_{\si \tau}(\la,\ga).
\]
For $\ga=\ga_0^\tau$ the proposition is now proved, since $E^-_\tau(x,\ga_0^\tau)=0$. Assume that $\ga=\ga_m^\tau$ for some $m \in \N$. For the anti-symmetric part we rewrite $E^-_\tau$ using Lemma \ref{lem:E-+} (note that $\de=\de_\tau$), and we use $\de(x) \De_\tau(x) =-(2x)^{-1} \De^+(x;u_0,t_0,t_1+1,u_1)$, then symmetrizing the integrand gives us
\[
\begin{split}
\langle& E^-_\tau(\cdot,\ga),\mathfrak E^-(\cdot,\la) \rangle_{\mathbf t^\tau}\\ 
=& (-1)^m\frac{\al_{\tau}\,\de_{\tau\si}(-\ga)\de_{\si}(\la)}{4\pi i} \int_{\mathcal C_\tau}  E^+(\cdot,\ga_{m-2}^\tau;t_0,u_0,t_1+1,u_1)\phi_\la (x;t_0,u_0,t_1+1,u_1) \\
& \mspace{160mu} \times \frac{\de(-x)-\de(x)}{2x} \De^+(x;u_0,t_0,t_1+1,u_1) dx \\
=& (-1)^{m+1} 2t_1\al_{\tau} \,\de_{\tau\si}(-\ga) \de_{\si}(\la) \big\langle E^+(\cdot,\ga_{m-2}^\tau;u_0,t_0,t_1+1,u_1), \phi_\la(\cdot;t_0,u_0,t_1+1,u_1) \big\rangle_{(u_0,t_0,t_1+1,u_1)}^+\\
=& (-1)^{m+1}2t_1 \al_{\tau}\,\de_{\tau\si }(-\ga)\de_{\si}(\la) \frac{ G_{\tau \si \tau}(\ga)}{G_{\tau \si \tau}(\ga_2^\tau)} E^+(\la,\ga_{m-2}^\tau;u_0,u_1,t_1+1,t_0)
\end{split}
\]
Here we used $\ga_{m-2}(u_0,t_0,t_1+1,u_1)=\pm(u_0+t_1+n)=\ga_{m}^\tau$ for $m=2n,2n-1$. From the explicit expressions for $\al_\tau$ and $G_{\tau \si \tau}(\ga_2^\tau)$ we find
\[
\frac{ \al_\tau}{G_{\tau \si \tau}(\ga_2^\tau)} = -\frac{ \al_{\si\tau}}{G_{\tau \si \tau}(\ga_0^\tau)}.
\]
Then adding the symmetric and the anti-symmetric parts we obtain 
\[
\begin{split}
\langle E_\tau&(\cdot,\ga), \mathfrak E(\cdot,\la) \rangle_{\mathbf t^\tau}\\
&= 2t_1 \frac{ G_{\tau \si \tau}(\ga) }{G_{\tau \si \tau}(\ga_0^\tau)} \Big(E^+_{\si \tau}(\la,\ga_{m}^\tau) + (-1)^m \al_{\si \tau} \,\de_{\si \tau \si}(-\ga)\de_{\si \tau}(\la) E^+(\la,\ga_{m-2}^\tau;u_0,u_1,t_1+1,t_0) \Big).
\end{split}
\]
Note that $\de_{\tau \si}= \de_{\si \tau \si}$ and $\de_\si = \de_{\si \tau}$, so applying Lemma \ref{lem:E-+} now gives \eqref{eq:id1}.
\end{proof}

Now we can express $\mathfrak F$, and therefore also $\mathcal F'$, as an integral transform. Our next goal is to write $\mathcal F'$ in the form $(\mathcal F'f)(\la) = \{ p, \mathcal E(\cdot,\la) \}$, where the bilinear form $\{\cdot,\cdot\}$ is such that the generators of $\mathcal H$ are symmetric with respect to this bilinear form.  The intertwining property $\mathcal F'\circ X= \si(X) \circ \mathcal F'$, $X \in \mathcal H$, then immediately gives transformation properties of the kernel $\mathcal E$ under the action of $\mathcal H$ and $\mathcal H_\si$.

Let $\Te$ be the weight function given by
\[
\Te(x) = G_\tau(x)^{-1} G(x)^{-1} \De(x) .
\]
By Euler's reflection formula for the $\Ga$-function we may write
\[
G_\tau(x)^{-1} G(x)^{-1} = \frac1{\pi^2} \sin \pi(d \pm x),
\]
and from this expression it is easy to see that $G_\tau^{-1} G^{-1}$ is $\mathcal W$-invariant. To the weight function $\Te$ we associate a bilinear form $\{\cdot,\cdot\}$ on $\mathcal AG_\tau$;
\[
\{f, g\}_{\mathbf t} = \frac{ 1}{2\pi i} \int_{\mathcal C} f(x) g(x) \Te(x) dx.
\]
\begin{lem} \label{lem:sym2}
The algebraic generators $z,T_0,T_1$ of $\mathcal H$ are symmetric with respect to $\{\cdot,\cdot\}_{\mathbf t}$.
\end{lem}
\begin{proof}
For $z$ this is trivial. For $i=0,1,$ we have $T_i(G_\tau(x)^{-1} G(x)^{-1} f) = G_\tau(x)^{-1} G(x)^{-1} T_if$, since $G_\tau(x)^{-1} G(x)^{-1}$ is $\mathcal W$-invariant. Since $T_i$ is symmetric with respect $\langle \cdot, \cdot \rangle_{\mathbf t}$ we have
\[
\{T_if, g\}_{\mathbf t} = \langle T_i f, G_\tau^{-1} G^{-1} g\rangle_{\mathbf t} = \langle f, T_i(G_\tau^{-1} G^{-1}g) \rangle_{\mathbf t} =\langle f, G_\tau^{-1} G^{-1} T_ig \rangle_{\mathbf t} = \{f, T_i g\}_{\mathbf t}. \qedhere
\]
\end{proof}
With the bilinear form $\{\cdot,\cdot\}$ we are going to construct a Fourier transform. The following function is the kernel in this Fourier transform.
\begin{Def} \label{def:nsWilson function}
We define the (non-symmetric) Wilson function by
\[
\mathcal E(x,\la) = G_\tau(x) G_{\si \tau}(\la) \mathfrak E(x,\la).
\]
\end{Def}
Observe that $\mathfrak E(x,\la)$ is an entire function in $x$ and $\la$, so the Wilson function is a meromorphic function in $x$ and $\la$ with simple poles coming from the Gaussians, located at $x= \pm (1-d+n)$, $\la=\pm(1-\tilde d +n)$, $n \in \Z_{\geq 0}$. With the non-symmetric Wilson function we now define an integral transform that is a Fourier transform associated to the duality isomorphism $\si$. Moreover, we can find the inverse and Plancherel-type formulas.
\begin{thm} \label{thm:Wil transform}
Let $\mathcal F$ be the linear operator defined by
\[
(\mathcal F f)(\la) = \{ f, \mathcal E(\cdot,\la) \}_{\mathbf t}, \qquad f \in \mathcal AG_\tau.
\]
Then

{\rm (a)} $\mathcal F$ maps $\mathcal AG_{\tau}$ into $\mathcal AG_{\si \tau}$.

{\rm (b)} $\mathcal F$ is a Fourier transform associated to $\si$.

{\rm (c)} $\mathcal F \circ \mathcal F_\si = (a+b)^2 \, id_{\mathcal AG_{\si \tau}}$ and $\mathcal F_\si \circ \mathcal F = (a+b)^2\, id_{\mathcal AG_\tau}$.

{\rm (d)} For $f_1,f_2 \in \mathcal AG_\tau$ and $g_1,g_2 \in \mathcal AG_{\si\tau}$,
\[
\{\mathcal Ff_1, \mathcal Ff_2\}_{\mathbf t^\si} = (a+b)^2\{f_1, f_2\}_{\mathbf t}, \qquad  \{\mathcal F_\si g_1, \mathcal F_\si g_2\}_{\mathbf t} = (a+b)^2\{g_1,g_2\}_{\mathbf t^\si}.
\]
\end{thm}
\begin{proof} 
{\rm (a)} We only need to check that the operator $\mathcal F$ is the same as $k\mathcal F'$, see \eqref{eq:def F}, for some non-zero constant $k$. Let $f \in \mathcal AG_\tau$. Observe that
\[
\De_\tau(x)=  G_\tau(x)  G(x)^{-1} \De(x) ,
\]
then using Proposition \ref{prop:kernel} we have
\[
\begin{split}
\mathcal K^{-1} \big(G_{\si\tau}\circ \mathfrak F \circ G_\tau^{-1}\big)f(\la) &= G_{\si\tau}(\la)\langle G_{\tau}^{-1}f, \mathfrak E(\cdot,\la) \rangle_{\mathbf t^\tau}\\
&= G_{\si\tau}(\la) \langle G_{\tau}^{-1}G^{-1}f, G_{\tau}\mathfrak E(\cdot,\la) \rangle_{\mathbf t}\\
&=\{ f ,\mathcal E(\cdot,\la) \}_{\mathbf t}.
\end{split}
\]
This is the desired result.

{\rm (b)} This follows from the proof of {\rm (a)} and Proposition \ref{prop:F}.

{\rm (c)} Let us denote $e_\ga(x)=G_\tau(x)E_\tau(x,\ga)$, for $\ga \in -\Ga_\tau$. From Proposition \ref{prop:kernel} we obtain
\begin{equation}  \label{eq:e1}
(\mathcal F e_\ga)(\la) =  (a+b) \frac{ G_{\tau \si \tau}(\ga) }{ G_{\tau \si \tau}(\ga_0^\tau)}G_{\si \tau}(\la) E_{\si \tau}(\la, \ga) =  (a+b) \frac{ G_{\tau \si \tau}(\ga) }{ G_{\tau \si \tau}(\ga_0^\tau)} e^\si_\ga(\la).
\end{equation}
Next we apply $\mathcal F_\si:\mathcal AG_{\si \tau} \rightarrow \mathcal A G_\tau$, which comes down to replacing $\mathbf t$ by $\mathbf t^\si$, then 
\[
\begin{split}
(\mathcal F_\si \circ \mathcal F)( e_\ga)(x)&= (a+b)\frac{ G_{\tau \si \tau}(\ga) }{ G_{\tau \si \tau}(\ga_0^\tau)}\,( \mathcal F_\si\, e_\ga^\si)(x) \\
&= (a+b)^2 \frac{G_{\tau \si \tau}(\ga) G_{\tau \si}(\ga) }{ G_{\tau \si \tau}(\ga_0^\tau) G_{\tau \si}(\ga_0^\tau)}\, e_\ga(x), 
\end{split}
\]
where we used $\mathbf t^{\si \tau \si \tau}= \mathbf t^{\tau \si}$. Writing out the Gaussians explicitly we have, for $m=2n$ or $m=2n-1$, 
\begin{equation} \label{eq:e2}
\frac{ G_{\tau \si}(\ga_{m}^\tau)}{ G_{\tau \si}(\ga_0^\tau)} = (-1)^n \frac{ (a+1-d)_n}{(b+c)_n}, \qquad \frac{ G_{\tau \si \tau}(\ga_{m}^\tau)}{ G_{\tau \si \tau}(\ga_0^\tau)} = (-1)^n \frac{(b+c)_n}{ (a+1-d)_n},
\end{equation}
so we have $(\mathcal F_\si \circ \mathcal F)(e_\ga)(x) =(a+b)^2 e_\ga(x)$. Since the non-symmetric Wilson polynomials form a basis for $\mathcal A$, the set $\{e_\ga\ |\ \ga \in -\Ga_\tau\}$ is a basis for $\mathcal AG_\tau$, and therefore $\mathcal F_\si \circ \mathcal F = (a+b)^2\, id_{\mathcal A G_\tau}$. In the same way we obtain $\mathcal F \circ \mathcal F_\si = (a+b)^2\, id_{\mathcal A G_{\si\tau}}$.

{\rm (d)} It is enough to show that the Plancherel formula is valid for the basis elements $e_\ga$ of $\mathcal AG_\tau$. By Theorem \ref{thm:FG} we have, for $\ga,\ga' \in -\Ga_\tau$,
\[
\{ e_\ga, e_{\ga'} \}_{\mathbf t} = \big\langle E_\tau(\cdot,\ga), E_\tau(\cdot,\ga') \big\rangle_{\mathbf t^\tau} 
= \de_{\ga,\ga'}\frac{ \langle 1, 1 \rangle_{\mathbf t^\tau} w_\tau(-\ga_0^\tau) }{w_\tau(-\ga)}.
\]
Writing out explicitly $\langle 1,1\rangle_{\mathbf t^{\tau}}$ in terms of $t_0,u_0,t_1,u_1$ using \eqref{eq:<1,1>} and \eqref{eq:abcd}, we find an expression that is invariant under $t_0 \leftrightarrow u_1$, so
\[
\langle 1,1 \rangle_{\mathbf t^{\tau}} = \langle 1,1 \rangle_{\mathbf t^{\si\tau}}.
\]
Using also the identity $G_{\tau \si \tau}(\ga)w_\tau(\ga)=G_{\tau \si}(\ga) w_{\si \tau}(\ga)$, see \eqref{eq:Gw}, we then obtain
\[
\{ e_\ga, e_{\ga'} \}_{\mathbf t} = \de_{\ga,\ga'}\frac{k_1(\ga) \langle 1, 1 \rangle_{\mathbf t^{\si\tau}} w_{\si\tau}(-\ga_0^\tau) }{w_{\si\tau}(-\ga)} = k_1(\ga) \big\langle E_{\si \tau}(\cdot,\ga), E_{\si \tau}(\cdot,\ga') \big\rangle_{\mathbf t^{\si \tau}} =k_1(\ga) \{e^\si_\ga,e^\si_{\ga'} \}_{\mathbf t^\si},
\]
where 
\[
k_1(\ga) = \frac{ G_{\tau \si}(\ga_0^\tau) G_{\tau \si \tau}(\ga)}{ G_{\tau \si}(\ga) G_{\tau \si \tau}(\ga_0^\tau)}.
\]
Now from \eqref{eq:e1} it follows that 
\[
\{ e_\ga, e_{\ga'} \}_{\mathbf t} = k_2(\ga) \{ \mathcal Fe_\ga, \mathcal Fe_{\ga'} \}_{\mathbf t^\si},
\]
where the factor $k_2(\ga)$ is given by
\[
k_2(\ga) = \frac{ k_1(\ga) G_{\tau \si \tau}(\ga_0^\tau)^2 }{ (a+b)^2 G_{\tau \si} \tau(\ga)^2} = (a+b)^{-2} \frac{G_{\tau \si}(\ga_0^\tau) G_{\tau \si \tau}(\ga_0^\tau)}{ G_{\tau \si}(\ga) G_{\tau \si \tau}(\ga) }.
\]
By \eqref{eq:e2} we have $k_2(\ga)=(a+b)^{-2}$, so we have proved the first Plancherel formula. The second formula follows from applying {\rm (c)}.
\end{proof}

Let us define the symmetric Wilson function by 
\[
\mathcal E^+(x,\la)=\big(C_+\mathcal E(\cdot,\la)\big)(x).
\]
Observe that $\mathcal E^+(x,\la)=G_\tau(x) G_{\si \tau}(\la)\phi_\la(x)$, with $\phi_\la$ defined by \eqref{lem:function phi}. So $\phi_\la$ can be considered as the analytic part of the symmetric Wilson function. With the symmetric kernel $\mathcal E^+$ we define the integral transform 
\[
(\mathcal F^+f)(\la) = \{f,\mathcal E^+(\cdot,\la)\}_{\mathbf t}^+, \qquad f \in \mathcal A_+G_\tau,
\]
where $\{\cdot,\cdot\}_{\mathbf t}^+$ is the bilinear form on $\mathcal A_+G_\tau$ defined by
\begin{gather*}
\{f,g\}_{\mathbf t}^+= \frac{1}{4 \pi  i} \int_{\mathcal C} f(x) g(x) \Te^+(x) dx,\\
\Te^+(x) = G_\tau(x)^{-1}G(x)^{-1} \De^+(x).
\end{gather*}
Then it not hard to verify that $\mathcal F^+:\mathcal A_+G_\tau \rightarrow \mathcal A_+G_{\si \tau}$ has the inverse $\mathcal F_\si^+$. Moreover, we have the Plancherel-type formula 
\[
\{\mathcal F^+f_1,\mathcal F^+f_2\}_{\mathbf t^\si}^+= \{f_1,f_2\}_{\mathbf t}^+, \qquad f_1,f_2 \in \mathcal A_+G_\tau,
\] 
so $\mathcal F^+$ is an isometry. 

Let us assume that the parameters $a,b,c,1-d$ are positive, or the non-real parameters occur in pairs of complex conjugates with positive real parts. Now the weights $\Te^+$ and $\De_\tau^+$ are positive on $i\R$, so $\{\cdot,\cdot\}_{\mathbf t}^+$ and $\langle \cdot,\cdot \rangle^+_{\mathbf t^\tau}$ define inner products. In general, if for some $\eps>0$ we have
\[
\int_\R e^{\eps|x|} d\mu(x)< \infty,
\]
then the moment problem for the measure $\mu$ is determinate, see \cite{Jeu03} and references therein, and therefore the polynomials are dense in the Hilbert space $L^2(\R, \mu)$. Using the asymptotic behavior \eqref{eq:asymp De} of $\De_\tau^+$, we find for $0<\eps<2\pi$
\[
\int_\R e^{\eps|x|} \De_\tau^+(ix) dx< \infty,
\]
so the polynomials are dense in $L^2(i\R,\De_\tau^+)$. In particular, the symmetric Wilson polynomials $E^+_\tau(x,\ga)$, $\ga \in \Ga^+$, form an orthogonal basis for $L^2(i\R, \De_\tau^+)^W$, the subspace of even functions in $L^2(i\R,\De_\tau^+)$. From this it follows that the set $\mathcal A_+G_\tau$ is dense in  $L^2(i\R,\Te^+)^W$. Therefore the integral transform $\mathcal F^+$ extends uniquely to a unitary operator $\tilde{\mathcal F}^+ : L^2(i\R,(4\pi i)^{-1}\Te^+)^W \rightarrow L^2(i\R,(4\pi i)^{-1}\Te^+_\si)^W$. Also, the operator $\mathfrak F^+ = G_{\si \tau}^{-1} \circ \mathcal F^+ \circ G_{\tau}$ extends uniquely to a unitary operator $\tilde{\mathfrak F}^+: L^2(i\R, (4\pi i)^{-1}\De_\tau^+)^W \rightarrow L^2(i\R, (4\pi i)^{-1}\De_{\si\tau}^+)^W$. This operator $\tilde{\mathfrak F}^+$ is precisely the Wilson function transform of type I, as defined in \cite{Gr03}.

\subsection{Properties of the Wilson function}
Let us finish this paper with a few nice properties of the Wilson function. From the property $\mathcal F\circ X=\si(X)\circ \mathcal F$ and Lemma \ref{lem:sym2} it follows immediately that the non-symmetric Wilson function satisfies
\begin{equation} \label{eq:X<->psiX}
\big(X \mathcal E(\cdot,\la) \big)(x) = \big( \psi(X) \mathcal E(x, \cdot)\big)(\la), \qquad X \in \mathcal H,
\end{equation}
where $\psi$ denotes the duality anti-isomorphism. In particular $\mathcal E(x,\la)$ is an eigenfunction of the Dunkl-Cherednik-type operators $Y$ and $Y^\si$ for eigenvalue $-\la$, respectively $-x$. From (the proof of) Proposition \ref{prop:A<->F} it follows that for $\la=\ga \in \Ga$ the Wilson polynomial $E(x,-\ga)$ also satisfies the transformation property \eqref{eq:X<->psiX}. This suggests that $\mathcal E(x,\ga) = k E(x,-\ga)$, $\ga \in \Ga$, for some non-zero constant $k$. Moreover, the dual Wilson function $\mathcal E_\si(\la,x)$ also satisfies \eqref{eq:X<->psiX}, which suggests that $\mathcal E(x,\la) = k \mathcal E_\si(\la,x)$ for some non-zero constant $k$.
\begin{thm} 
The Wilson function $\mathcal E(x,\la)$ satisfies:\\
{\rm (a)} $\displaystyle
\mathcal E(x,\ga) = \frac{\Ga(1-a-d)}{\Ga(a+b) \Ga(a+c)} E(x,-\ga),$ for $\ga \in \Ga$\\
{\rm (b)} $\mathcal E(x,\la) = \mathcal E_\si(\la,x)$
\end{thm} 
So $\mathcal E(x,\la)$ is a meromorphic continuation of the non-symmetric Wilson polynomial in its degree, and like the Wilson polynomial, $\mathcal E(x,\la)$ also satisfies a duality property.
\begin{proof}
Consider first the function $\phi_\la(x)$ defined by \eqref{def:phi}. Using $a+e=\tilde a +\tilde e$ for $e=b,c,d$ and $e-d=\tilde e-\tilde d$ for $e=b,c$ the duality property for $\phi_\la$ follows directly from \eqref{def:phi}:
\[
\phi_\la(x;\mathbf t) = \phi_x(\la;\mathbf t^\si).
\]
Now the duality property for the Wilson function follows from Definition \ref{def:nsWilson function} and the definition of the kernel $\mathfrak E$ in Proposition \ref{prop:kernel}.

For the reduction to the Wilson polynomials observe the function $\Ga(\tilde a \pm \la)^{-1}$ has zeros at $\la = \pm(\tilde a+n)$, $n \in \Z_{\geq 0}$. Thefore the second $_4F_3$-function in \eqref{def:phi} vanishes for $\la \in \Ga$, and then we see that
\[
\phi_{\ga}(x) = \frac{\Ga(1-a-d)}{\Ga(a+b) \Ga(a+c) G_\tau(x) G_{\si\tau}(\ga)} E^+(x,\ga), \qquad \ga \in \Ga.
\]
Now (\rm a) follows from Definition \ref{def:nsWilson function} and Lemma \ref{lem:E-+}.
\end{proof}


\begin{thebibliography}{99}
\bibitem{AW85} R. Askey, J. Wilson, \emph{Some basic hypergeometric orthogonal polynomials that generalize Jacobi polynomials},  Mem. Amer. Math. Soc. \textbf{54} (1985).

\bibitem{Ch95} I. Cherednik, \emph{Nonsymmetric Macdonald polynomials},   Internat. Math. Res. Notices \textbf{10} (1995), 483-515.

\bibitem{Ch97} I. Cherednik, \emph{Difference Macdonald-Mehta conjecture}, Internat. Math. Res. Notices  \textbf{10} (1997), 449-467.

\bibitem{Ch97A} I. Cherednik, \emph{Inverse Harish-Chandra transform and difference operators}, Internat. Math. Res. Notices \textbf{15} (1997), 733-750.

\bibitem{CM02} I. Cherednik, Y. Markov, \emph{Hankel transform via double Hecke algebra}, in: Iwahori-Hecke algebras and their representation theory (Martina-Franca, 1999),  1-25, Lecture Notes in Math., \textbf{1804}, Springer, Berlin, 2002.

\bibitem{Gr03} W. Groenevelt, \emph{The Wilson function transform}, Internat. Math. Res. Notices. \textbf{52} (2003), 2779-2817.

\bibitem{ILVW90} M.E.H. Ismail, J. Letessier, G. Valent, J. Wimp, \emph{Two families of associated Wilson polynomials}, Canad. J. Math. \textbf{42} (1990),  659-695.

\bibitem{Jeu03} M. de Jeu, \emph{Determinate multidimensional measures, the extended Carleman theorem and quasi-analytic weights},  Ann. Probab. \textbf{31} (2003), 1205-1227.

\bibitem{KS01} E. Koelink, J.V. Stokman,\emph{The Askey-Wilson function transform},  Internat. Math. Res. Notices \textbf{22} (2001), 1203-1227.

\bibitem{Ko85} T.H. Koornwinder, \emph{Special orthogonal polynomial systems mapped onto each other by the Fourier-Jacobi transform}, in: Orthogonal Polynomials and Applications (Bar-le-Duc, 1984), 174-183, Lecture Notes in Math., 1171, Springer, Berlin, 1985.

\bibitem{Ko92}  T.H. Koornwinder, \emph{Askey-Wilson polynomials for root systems of type $BC$}, Contemp. Math., \textbf{138} (1992), 189-204.

\bibitem{Ma91}  D.R. Masson, \emph{Associated Wilson polynomials}, Constr. Approx. \textbf{7} (1991), 521-534. 

\bibitem{Mac03} I.G. Macdonald, \emph{Affine Hecke Algebras and Orthogonal Polynomials}. Cambridge Tracts in Mathematics, 157. Cambridge University Press,  2003.

\bibitem{No95} M. Noumi, \emph{Macdonald-Koornwinder polynomials and affine Hecke rings}, Surikaisekikenkyusho Kokyuroku  \textbf{919} (1995), 44-55, (in Japanese).

\bibitem{NS04} M. Noumi, J.V. Stokman, \emph{Askey-Wilson polynomials: an affine Hecke algebra approach}, in: Laredo Lectures on Orthogonal Polynomials and Special Functions, 111-144, Adv. Theory Spec. Funct. Orthogonal Polynomials, Nova Sci. Publ., Hauppauge, NY, 2004.

\bibitem{Op95}  E.M. Opdam, \emph{Harmonic analysis for certain representations of graded Hecke algebras}, Acta Math. \textbf{175} (1995), 75-121.

\bibitem{Sa99} S. Sahi, \emph{Nonsymmetric Koornwinder polynomials and duality},  Ann. of Math. \textbf{150} (1999), 267-282.

\bibitem{St00} J.V. Stokman, \emph{Koornwinder polynomials and affine Hecke algebras}, Internat. Math. Res. Notices \textbf{19} (2000), 1005-1042.

\bibitem{St03} J.V. Stokman, \emph{Difference Fourier transforms for nonreduced root systems}, Selecta Math. (N.S.)  \textbf{9} (2003), 409-494.

\bibitem{Wi80} J.A. Wilson, \emph{Some hypergeometric orthogonal polynomials},  SIAM J. Math. Anal. \textbf{11} (1980), 690-701.

\bibitem{Zh05} G. Zhang, \emph{Spherical transform and Jacobi polynomials on root systems of type BC}, math.RT/0503735.
\end{thebibliography}
\end{document}